\UseRawInputEncoding
\documentclass{article}

\usepackage{color,graphicx,times}

\makeatother

\begin{document}
 
\date{}

\title{\Large\bf Centrality and the KRH Invariant } 

\author{Louis H. Kauffman \\
and\\
David Radford\\
Department of Mathematics, Statistics and Computer Science \\
University of Illinois at Chicago \\
851 South Morgan Street\\
Chicago, IL 60607-7045\\
and\\
Stephen Sawin\\
Mathematics Department\\
Fairfield University\\
Fairfield, Connecticut 06430-5195\\
}

\maketitle

\thispagestyle{empty}

\section{Introduction}

\noindent{\bf Remark.} The body of this paper is identical with \cite{Centrality} published in 1998. In this document, the references have
been updated and applications to virtual knot theory have been added in Section 8. The formalism of the invariants considered here applies directly to rotational virtual knot theory
as we have formulated it in \cite{RV}.\\

The purpose of this paper is to discuss some consequences of a method of
defining invariants of links and 3-manifolds intrinsically in terms of right integrals on certain Hopf algebras. We call such an invariant of 3-manifolds a Hennings invariant \cite{Hennings}. The work reported in this paper has as its background \cite{KR}, \cite{KandR}, \cite{K-Hopf}.
\vspace{3mm}

Hennings invariants were originally defined using oriented links. It is not necessary to use invariants that are dependent on link orientation to define 3-manifold invariants via surgery and Kirby calculus. The invariants discussed in this paper are formulated for unoriented links. This results in a simplification and conceptual clarification of the relationship of Hopf algebras and the link invariants defined by Hennings.
\vspace{3mm}

We show in \cite{KandR} that invariants defined in terms of right integrals, as considered in this paper, are distinct from the invariants of Reshetikhin and Turaev \cite{RT}. We show that the Hennings invariant is non-trivial for
the quantum group $U_{q}(sl_{2})'$ when $q$ is an fourth root of unity. The Reshetikhin Turaev invariant is trivial at this quantum group and root of unity. The Hennings invariant distinguishes all the Lens spaces $L(n,1)$ from one another at this root of unity. This proves that there is non-trivial topological information in the non-semisimplicity of $U_{q}(sl_{2})'$. This non-triviality result has also been obtained by Ohtsuki \cite{Ohtsuki}.
\vspace{3mm}

The reader interested in comparing the approach of this paper with other ways to look at quantum link invariants will enjoy looking at the references \cite{KandP}, \cite{Kerler}, \cite{KirbyMelvin}, \cite{Lawrence}, \cite{Luba}, \cite{Radford-trace}, \cite{Radford-gen}, \cite{Resh}, \cite{RT}, \cite{Witten}. In particular, the method we use to write link invariants directly in relation to a Hopf algebra is an analog of the construction in \cite{Lawrence} and it is a generalization of the formalism of \cite{Resh} and \cite{RT}. The papers \cite{Kerler}, \cite{Luba}, \cite{Kuperberg1}, \cite{Kuperberg2} consider
categorical frameworks that also use right integrals on Hopf algebras. More information about Hopf algebras in relation to our constructions can be found in \cite{Radford-trace} and \cite{Radford-gen}. The book \cite{KandP} contains background material on link invariants from many points of view, including a sketch of the method taken in \cite{Witten}. It is an open question whether there is a natural quantum field theoretic interpretation of the three-manifold invariants discussed in this paper. \vspace{3mm}

The present paper is a review of the structure of these invariants and it emphasizes how the
framework developed for this study of invariants can be regarded as a construction of a natural category associated with a (finite dimensional) quasitriangular Hopf algebra. The construction provides a functor from the category of tangles to this category associated with the Hopf algebra. In this context we obtain an elegant proof that the image under this functor of 1-1 tangles gives elements in the center of the Hopf algebra. This constitutes a non-trivial application of these categories to the structure of Hopf algebras. \vspace{3mm}

The paper is organized as follows. Section 2 recalls Hopf algebras, quasitriangular Hopf algebras and ribbon Hopf algebras. Section 3 discusses the conceptual setting of the invariant via the different categories of tangles, immersions and morphisms associated with the Hopf algebra. In Section 4 we discuss the diagrammatic and categorical structures associated with traces and integrals on Hopf algebras. In Section 5 we show that traces of the kind discussed in Section 3 can be constructed from right integrals, and that these traces yield invariants of the 3-manifolds obtained by surgery on the links. Section 6 details the application to centrality, giving a very simple proof that the elements $F(T)$ of the Hopf algebra that are images of $1-1$ tangles under our functor from tangles to Hopf algebras are
in the center of the Hopf algebra. Section 7 points out that the centrality proof in Section 6 is actually constructive at the algebra level, giving specific proofs of centrality for each example. Furthermore, a direct analysis of the relationship of the combinatorics and
the algebra reveals another proof of centrality based on pushing algebraic
beads around the diagram. We describe
and prove this result, raising the question at the end of a full characterisation of central elements in the Hopf algebra. \vspace{3mm}

\noindent
{\bf Acknowledgment.} Louis Kauffman thanks the National Science Foundation for support of this research under grant number DMS-9205277. David Radford thanks the National Science Foundation for support under grant number DMS-870-1085. Steve Sawin thanks
the National Science Foundation for support under NSF Postdoctoral Fellowship number 23068.
\vspace{3mm}

\section{Algebra}
Recall that a Hopf algebra $A$ \cite{Sweedler} is a bialgebra over a commutative ring $k$ that has an associative multiplication and a coassociative comultiplication and is equipped with a counit, a unit
and an antipode. The ring $k$ is usually taken to be a field. $A$ is an algebra with multiplication $m:A \otimes A \longrightarrow A$.
The associative law for m is expressed by the equation $m(m \otimes1_{A}) = m(1_{A} \otimes m)$ where $1_{A}$ denotes the identity map on A.
\vspace{3mm}

The coproduct
$\Delta :A \longrightarrow A \otimes A$ is an algebra homomorphism and is
coassociative in the sense that
$(\Delta \otimes 1_{A})\Delta = (1_{A} \otimes \Delta) \Delta.$ \vspace{3mm}

The unit is a mapping from $k$ to $A$ taking $1$ in $k$ to $1$ in $A$, and thereby defining an action of $k$ on $A.$ It will be convenient to just identify the units in $k$ and in $A$, and to ignore
the name of the map that gives the unit. \vspace{3mm}

The counit is an algebra mapping from $A$ to $k$ denoted by $\epsilon :A \longrightarrow k.$ The following formulas for the counit dualize the structure inherent in the unit:
$(\epsilon \otimes 1_{A}) \Delta = 1_{A} = (1_{A} \otimes \epsilon) \Delta.$
\vspace{3mm}

It is convenient to write formally
$$\Delta (x) = \sum x_{1} \otimes x_{2} \in A \otimes A$$
to indicate the decomposition of the coproduct of $x$ into a sum of first and second factors in the two-fold tensor product of $A$ with itself. We shall often drop the summation sign and write
$$\Delta (x) = x_{1} \otimes x_{2}.$$
\vspace{3mm}

The antipode is a mapping $s:A \longrightarrow A$ satisfying the equations
$m(1_{A} \otimes s) \Delta (x) = \epsilon (x)1,$ and $m(s \otimes 1_{A}) \Delta (x)= \epsilon (x)1$ where 1 on the right hand side of these equations denotes the unit of $k$ as identified with the unit of $A.$ It is a consequence of this definition that $s(xy) = s(y)s(x)$ for all $x$ and $y$ in A. \vspace{3mm}

A quasitriangular Hopf algebra $A$ \cite{Drinfeld} is a Hopf algebra with an element $\rho \in A \otimes A$ satisfying
the following equations:
\vspace{3mm}

\noindent
1) $\rho \Delta = \Delta' \rho$ where $\Delta'$ is the composition
of $\Delta$ with the map on
$A \otimes A$ that switches the two factors. \vspace{3mm}

\noindent
2) $$\rho_{13} \rho_{12} = (1_{A} \otimes \Delta) \rho,$$ $$\rho_{13} \rho_{23} = (\Delta \otimes 1_{A})\rho.$$ \vspace{3mm}

\noindent
{\bf Remark.} The symbol $\rho_{ij}$ denotes the placement of the first and second tensor factors of $\rho$ in the $i$ and
$j$ places in a triple tensor product. For example, if $\rho = \sum e \otimes e'$ then $$\rho_{13} = \sum e \otimes 1_{A} \otimes e'.$$
\vspace{3mm}

These conditions imply that $\rho$ has an inverse, and that

$$ \rho^{-1} = (1_{A} \otimes s^{-1}) \rho = (s \otimes 1_{A}) \rho .$$

It follows easily from the axioms of the quasitriangular Hopf algebra that $\rho$ satisfies the Yang-Baxter equation

$$\rho_{12} \rho_{13} \rho_{23} = \rho_{23} \rho_{13} \rho_{12}.$$

A less obvious fact about quasitriangular Hopf algebras is that there exists an element $u$ such that $u$ is invertible and $s^{2}(x) = uxu^{-1}$ for all $x$ in $A.$ In fact, we may take $u = \sum s(e')e$ where $\rho = \sum e \otimes e'.$ As we shall see, this result, originally due to Drinfeld \cite{Drinfeld}, follows from the diagrammatic categorical context of this paper. \vspace{3mm}

An element $G$ in a Hopf algebra is said to be {\em grouplike} if $\Delta (G) = G \otimes G$ and $\epsilon (G)=1$ (from which it follows that $G$ is invertible and $s(G) = G^{-1}$). A quasitriangular Hopf algebra is said to be a {\em ribbon Hopf algebra} \cite{RTG}, \cite{KR} if there exists a grouplike element $G$ such that (with $u$ as in the previous paragraph) $v = G^{-1}u$ is in
the center of $A$ and $s(u) = G^{-1}uG^{-1}$. We call G a special grouplike element of A.
\vspace{3mm}

Since $v=G^{-1}u$ is central, $vx=xv$ for all x in A. Therefore $G^{-1}ux = xG^{-1}u.$ We know that $s^{2}(x) = uxu^{-1}.$ Thus $s^{2}(x) =GxG^{-1}$ for all $x$ in $A.$ Similarly, $s(v) = s(G^{-1}u) = s(u)s(G^{-1})=G^{-1}uG^{-1}G =G^{-1}u=v.$ Thus the square of the
antipode is represented as conjugation by the special grouplike element in a ribbon Hopf algebra, and the central element $v=G^{-1}u$ is invariant under the antipode.
\vspace{3mm}

\section{Tangle Categories and Hopf Algebra Categories} We now describe the categories that are the contexts for the results of this paper. These categories span the gamut from topology to algebra. At one end we have the category of tangles. At the
other end we have a natural category associated to any algebra, where the elements of the algebra become the morphisms of the category. This section will describe the categories that we need and the relevant functors between them.
\vspace{3mm}

\subsection{The Tangle Category}
We begin with the (unoriented) tangle category, which we refer to as {\em
Tang}. This category is the main topological category that we use. It encompasses knots, links, braids and their generalizations, known as {\em tangles}. All of the usual objects from the point of view of a topologist become morphisms in $Tang.$ The objects in $Tang$ consist of formal finite tensor products of a basic object $V$ with itself and with another object $k$. These tensor products are all associative, and obey the following rules: $V \otimes k$ = $k \otimes V$ = $V$
and $k \otimes k$ = $k$. Of course $V$ is a formal analogue of a module over a ring $k$.
\vspace{3mm}

While the objects in the tangle category are very simple, the morphisms are quite complex. Each morphism in the tangle category consists in a link diagram with free ends which is transverse with respect to
a given direction in the plane. (This special direction will be called the {\em vertical} direction.) The transversality of the diagram to this vertical direction means that any given line perpendicular to the vertical direction intersects the diagram either tangentially at a maximum or a minimum, at non-zero angle for any other strand. We shall further assume that any given perpendicular intersects the diagram at at most one crossing. With these stipulations the free ends of the diagram occur at either its top or its bottom (top and bottom taken with respect to the designated vertical direction).We shall assume that all the top ends occur along the same perpendicular, and that all the bottom ends occur along another perpendicular to the vertical. To each of these two rows of diagram ends is assigned a tensor product of copies of $V$, one for each end. In the tangle category, the diagram is a morphism from the lower tensor product to the upper tensor product. Thus a diagram with $n$ lower ends and $m$ upper ends is a morphism from $V ^{\otimes ^{n}}$ to $V^{ \otimes ^{m}}$. If the top of a diagram has no free ends, then its range is $k$. If the bottom of a diagram has no free ends, then its domain is $k$. A diagram is said to be {\em closed} if it has no free ends. Thus a closed diagram is a morphism from $k$ to $k$. If $A$ and $B$ are morphisms in the tangle category with $Range(B) = Domain(A)$, then the composition of $A$ and $B$ is denoted $AB$. (The reader should note that we have taken this left-right convention for the composition of morphisms in the tangle category. The left-right convention is opposite to that usually adopted for function composition. In composing functions we use the usual convention; $gf(x) = g(f(x))$ when $Range(f) = Domain(g)$.
\vspace{3mm}

Let $U$ and $V$ be morphisms in $Tang.$ We define their {\em tensor product}, $U \otimes V$, to be the tangle obtained from the tangles $U$ and $V$ by juxtaposing them disjointly side by side, with $U$ to the left of $V$. In other words, the inputs to $U \otimes V$ consist in the inputs to $U$ followed by the inputs to $V$, and similarly for the outputs.) The domain of $U \otimes V$ is the tensor product of the domains of $U$ and $V$, and the range is the tensor product of the ranges. This makes $Tang$ into a tensor category. Note that a tangle consisting in a single upward-moving line is taken to be the identity map from $V$ to $V$, and hence a tangle consisting in $n$ parallel lines is the identity map on the $n$-fold tensor product of $V$ with itself. \vspace{3mm}

It is not hard to see that every morphism in the tangle category is a composition of the elementary morphisms $Cup$,$Cap$, $R$ and $L.$
$Cup$ denotes a tangle with
no inputs and two outputs that are connected by a single arc that forms a minimum. $Cap$ denotes a tangle with no outputs and two inputs that are connected by a single arc that forms a maximum. $R$ denotes a crossing of two arcs so that the overcrossing line goes upward from right to left. $L$ denotes a crossing of two arcs so that the overcrossing line goes upward from left to right. Both $R$ and $L$ are tangles with two inputs and two outputs. Thus, as morphisms in $Tang$, we have $$Cup: k \longrightarrow V \otimes V$$ and
$$Cap: V \otimes V \longrightarrow k.$$ $R$ and $L$ are each morphisms of the form $V \otimes V
\longrightarrow V \otimes V.$ See Figure 1. In the next paragraphs, we will discuss the axioms that will be imposed on these generating morphisms. \vspace{3mm}

\begin{figure}[htb]
     \begin{center}
     \begin{tabular}{c}
     \includegraphics[width=7cm]{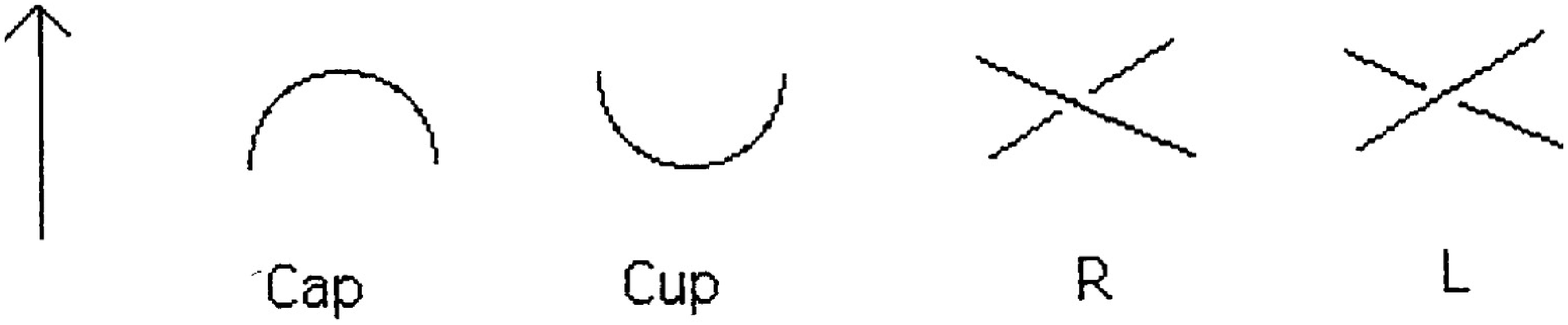}
     \end{tabular}
     \caption{\bf Cups, Caps and Crossings}
     \label{F1}
\end{center}
\end{figure}

Now we discuss the equivalence relation on the morphisms in $Tang$. This equivalence corresponds directly to regular isotopy of link diagrams and tangles arranged with respect to a ``vertical" direction. For this reason, we shall discuss this equivalence relation first in topological terms, and then transfer the discussion to the category.
\vspace{3mm}

As described above, any link diagram or tangle can be arranged to be transversal to a given direction (designated as vertical) in the plane. We shall call such diagrams {\em vertical} diagrams. Once we assume that the diagrams are so given, it is necessary to add two more moves to the classical list of Reidemeister moves \cite{Reidemeister} in order to insure that intermediate stages in an isotopy remain vertical. The resulting four vertical moves (we do not use the classical first Reidemeister move, since the intent is to model regular isotopy) are illustrated in Figure 2.
\vspace{3mm}

\begin{figure}[htb]
     \begin{center}
     \begin{tabular}{c}
     \includegraphics[width=7cm]{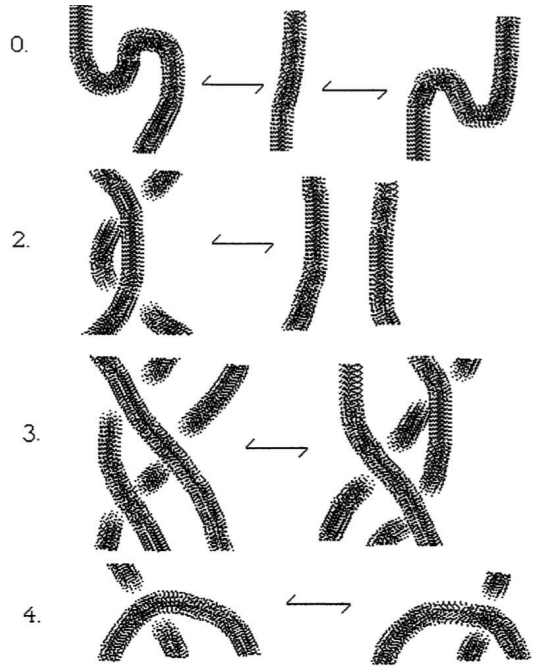}
     \end{tabular}
     \caption{\bf Reidemeister Moves With Respect To A Vertical Direction}
     \label{F2}
\end{center}
\end{figure}

\vspace{3mm}

Move $0$ comprises the cancellation of adjacent maxima and minima.
Move $2$ can be regarded as the cancellation of crossings of opposite type.
Move $3$ is the basic braiding identity. Move $4$ is a ``switchback" move that exchanges a crossing next to a maximum (minimum) for the maximum (minimum) next to the opposite crossing. Each of the moves can be regarded as a relation on the generating morphisms $Cup$, $Cap$, $L$ and $R$. Specifically, here are the corresponding algebraic statements of these moves: \vspace{3mm}

\noindent
0. $(Cup \otimes 1_{V})(1_{V} \otimes Cap) = 1_{V} \otimes 1$, 
$(1_{V} \otimes Cup)(Cap \otimes 1_{V}) =1_{V}.$ \\ \vspace{2mm}

\noindent
2. $RL = LR = 1_{V} \otimes 1_{V}.$\\
\vspace{2mm}

\noindent
3. $(R \otimes 1_{V})(1_{V} \otimes R)(R \otimes 1_{V}) = (1_{V} \otimes R)(R \otimes 1_{V})(1_{V} \otimes R).$\\ \vspace{2mm}

\noindent
4. $(L \otimes 1_{V})( 1_{V} \otimes Cap) = (1_{V} \otimes R)(Cap \otimes 1_{V}).$\\
\vspace{3mm}

Each of these equations is taken in the context of the identifications $k \otimes k = k$ and $V \otimes k = V = k \otimes V$. The notation $1_{V}$ stands for the identity map on $V$. This completes the description of the category $Tang.$
\vspace{3mm}

In discussing $Tang$ we shall continue to use the topological terminology {\em tangle} for a morphism in the category. An {\em $n-m$ tangle} is a tangle with $n$ inputs and $m$ outputs. Thus a $1-1$ tangle is any map from $V$ to $V$ in $Tang$, and a $0-0$ tangle is any knot or link arranged with respect to the vertical to give a morphism from $k$ to $k$.
If $T$ is an $n-m$ tangle in $Tang$ we define $\Delta (T)$ to be the $2n-2m$ tangle obtained by replacing every strand of $T$ by two parallel copies of that strand. We leave it as an exercise for
the reader to express $\Delta$ on the generating morphisms. As we shall see, $\Delta$ is an analog of the coproduct in a Hopf algebra. There is also an analog of the antipode in a Hopf algebra, defined in $Tang$ and taking an $n-m$ tangle $T$ to an $m-n$ tangle $S(T)$. The functor $S$ is defined by adding caps on the left and cups on the right of the tangle $T$ so that the inputs and outputs are reversed. See Figure 3 for an illustration of the tangle antipode $S$. In $Tang$ we have $S^{2} = I$ where $I$ denotes the identity functor on tangles. This is a precursor to the special nature of the elements of the Hopf algebra category that will be the image of our functor from the tangle category to the Hopf algebra category. \vspace{3mm}

\begin{figure}[htb]
     \begin{center}
     \begin{tabular}{c}
     \includegraphics[width=7cm]{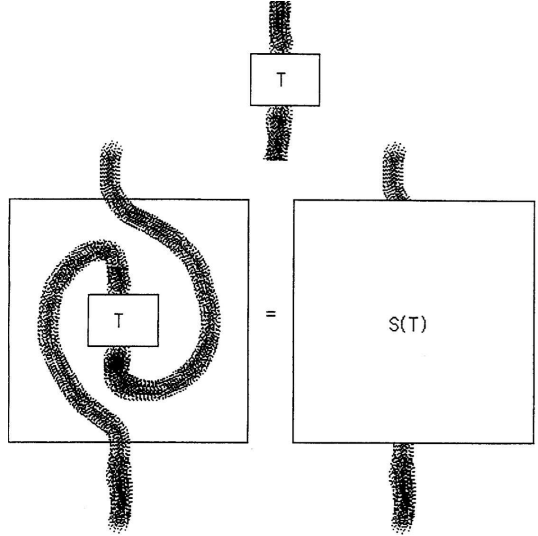}
     \end{tabular}
     \caption{\bf  Tangle Antipode}
     \label{F3}
\end{center}
\end{figure}

\subsection{The Immersion Category}
The {\em Immersion Category}, denoted {\em Flat}, is a quotient of the tangle category where we identify the maps $L$ and $R$ with each other.
In this category, let $P$ denote the equivalence class of $L=R$. The morphisms in $Flat$ are represented by tangle diagrams that are immersed in the plane (that is each curve in the diagram is the image of an immersion and distinct curves intersect transversely). The axioms for equivalence of morphisms in $Flat$ correspond to regular homotopy of flat tangles. The Whitney-Graustein Theorem \cite{Whitney} applies to this category. The Whitney-Graustein Theorem states that any immersed curve in the plane is regularly homotopic to a simple closed curve that is decorated with a string of curls, as indicated in Figure 4. A {\em curl} is a $1-1$ tangle in $Flat$ of the form $G=(1_{V} \otimes Cup)P(1_{V} \otimes Cap)$ or $G^{-1} = (Cup \otimes 1_{V})P(Cap \otimes 1_{V})$. The second curl is denoted by $G^{-1}$ because the equations $GG^{- 1} = 1_{V} = G^{-1}G$ hold in $Flat$ as illustrated in Figure 5. These equations are the categorical analog of the so- called "Whitney trick". Whitney's Theorem tells us that any $1-1$ flat tangle with a single strand is equivalent to an integer power of $G$. The exponent is the Whitney degree of this curve, oriented from input to output. (The Whitney degree is the total turn of the tangent vector to the oriented curve. If the curve is a $1-1$ tangle then the Whitney degree of the identity tangle is zero. If the curve is a single strand $0-0$ tangle, then the Whitney degree depends upon the choice of orientation, with a clockwise oriented circle giving degree $1$. \vspace{3mm}

\begin{figure}[htb]
     \begin{center}
     \begin{tabular}{c}
     \includegraphics[width=7cm]{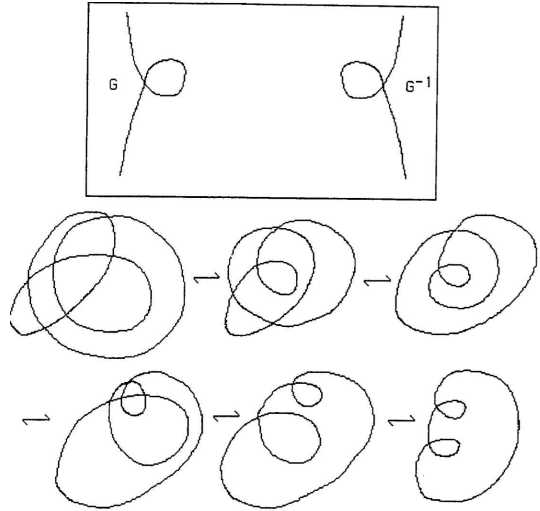}
     \end{tabular}
     \caption{\bf  Immersed Curves}
     \label{F4}
\end{center}
\end{figure}

 \begin{figure}[htb]
     \begin{center}
     \begin{tabular}{c}
     \includegraphics[width=7cm]{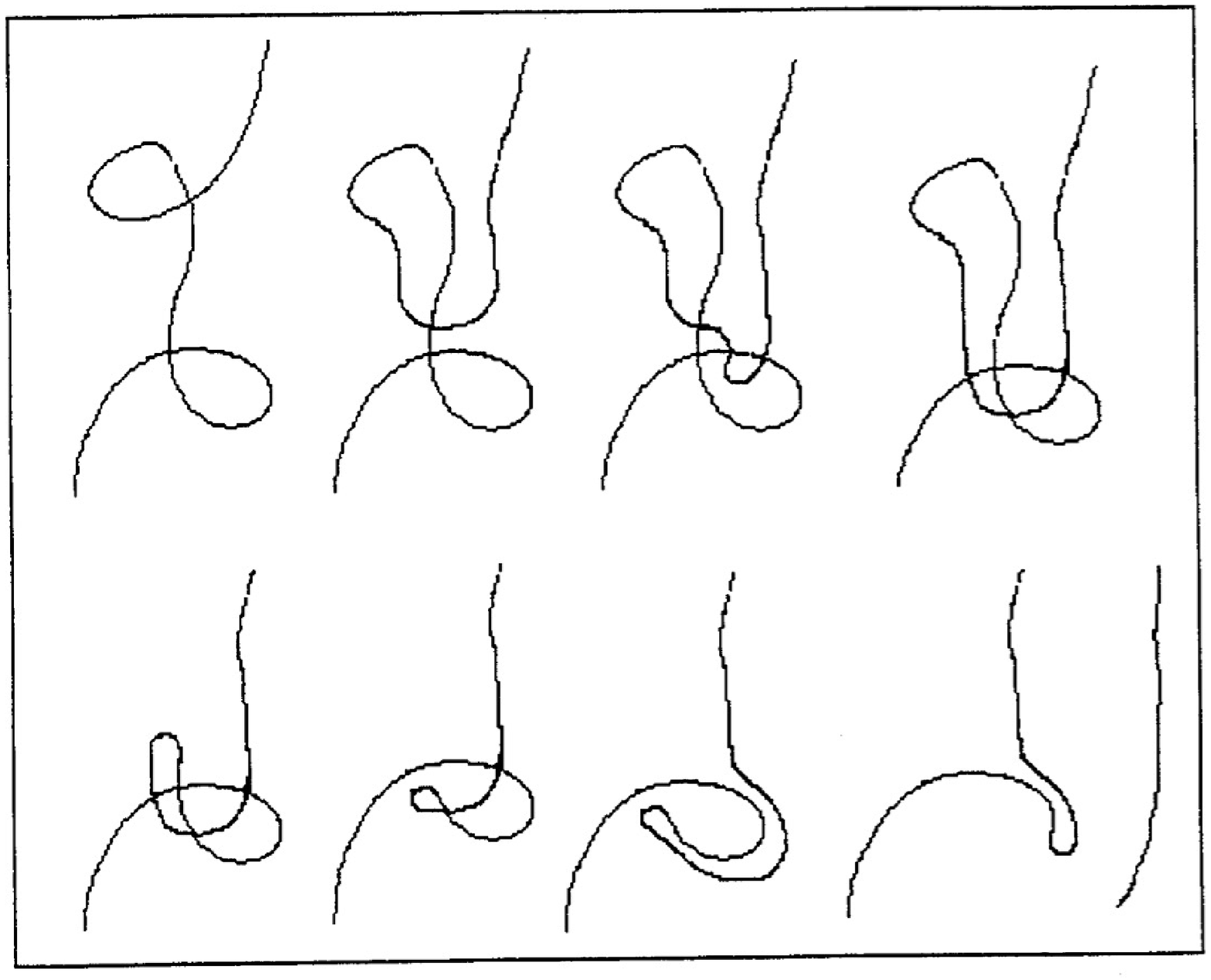}
     \end{tabular}
     \caption{\bf Whitney Trick}
     \label{F5}
\end{center}
\end{figure}

Note that in $Flat$ the morphism
$P:V \otimes V \longrightarrow V \otimes V$ has the formal properties of a permutation of the factors of $V \otimes V$. For example $P^{2} = 1_{V} \otimes 1_{V}$. In computing the categorical antipode $S$ on the elements $G$ and $G^{-1}$, we find that $\Delta(G) = G \otimes G$ and $\Delta(G^{-1}) = G^{-1} \otimes G^{-1}$, as shown in Figure 6. This means that $G$ behaves as a formal group-like element in the category $Flat$. Note also that $S(G) = G^{-1}$, as shown in Figure 6.
\vspace{3mm}

 \begin{figure}[htb]
     \begin{center}
     \begin{tabular}{c}
     \includegraphics[width=7cm]{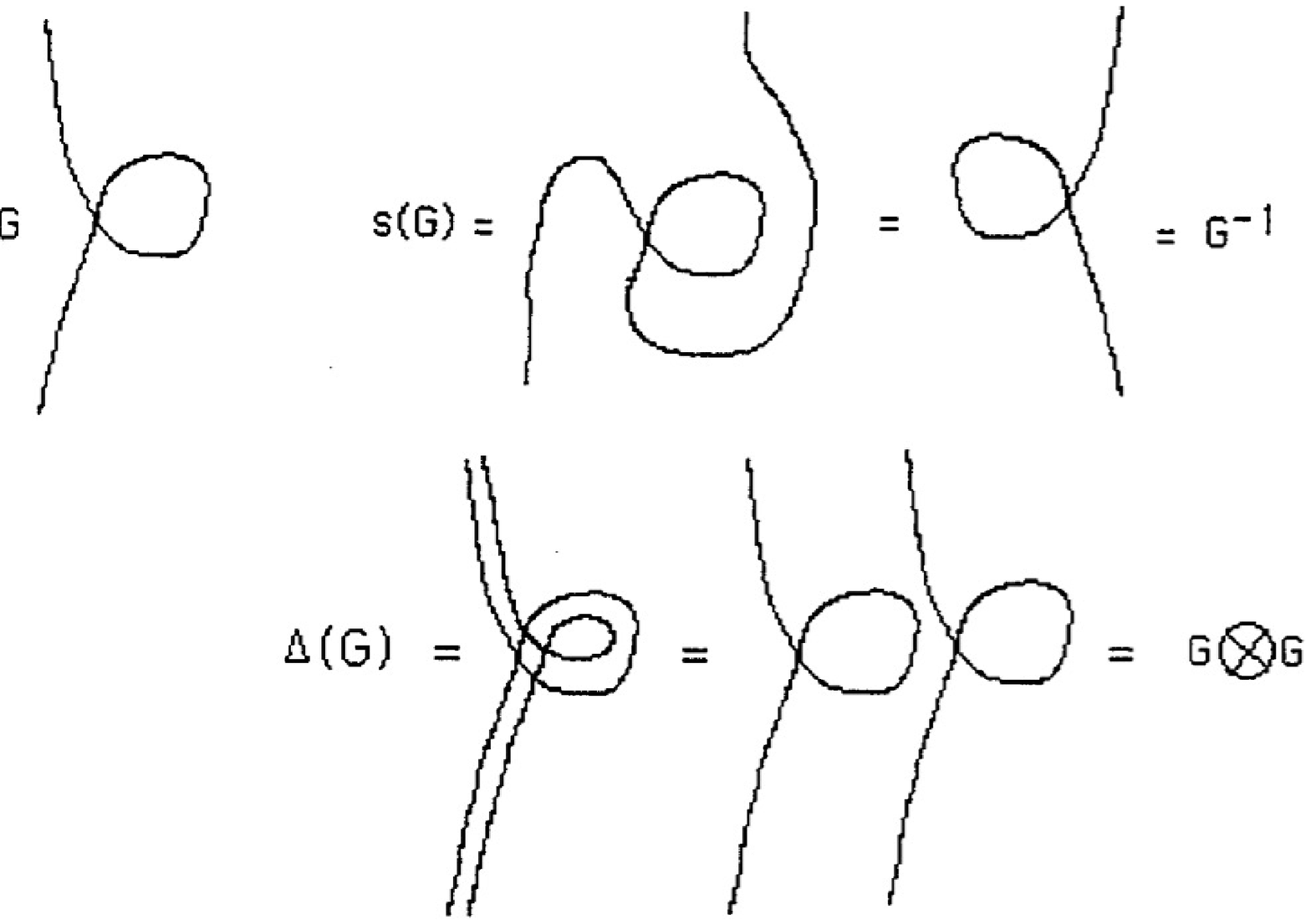}
     \end{tabular}
     \caption{\bf G is a formal grouplike element.}
     \label{F6}
\end{center}
\end{figure}

 \begin{figure}[htb]
     \begin{center}
     \begin{tabular}{c}
     \includegraphics[width=7cm]{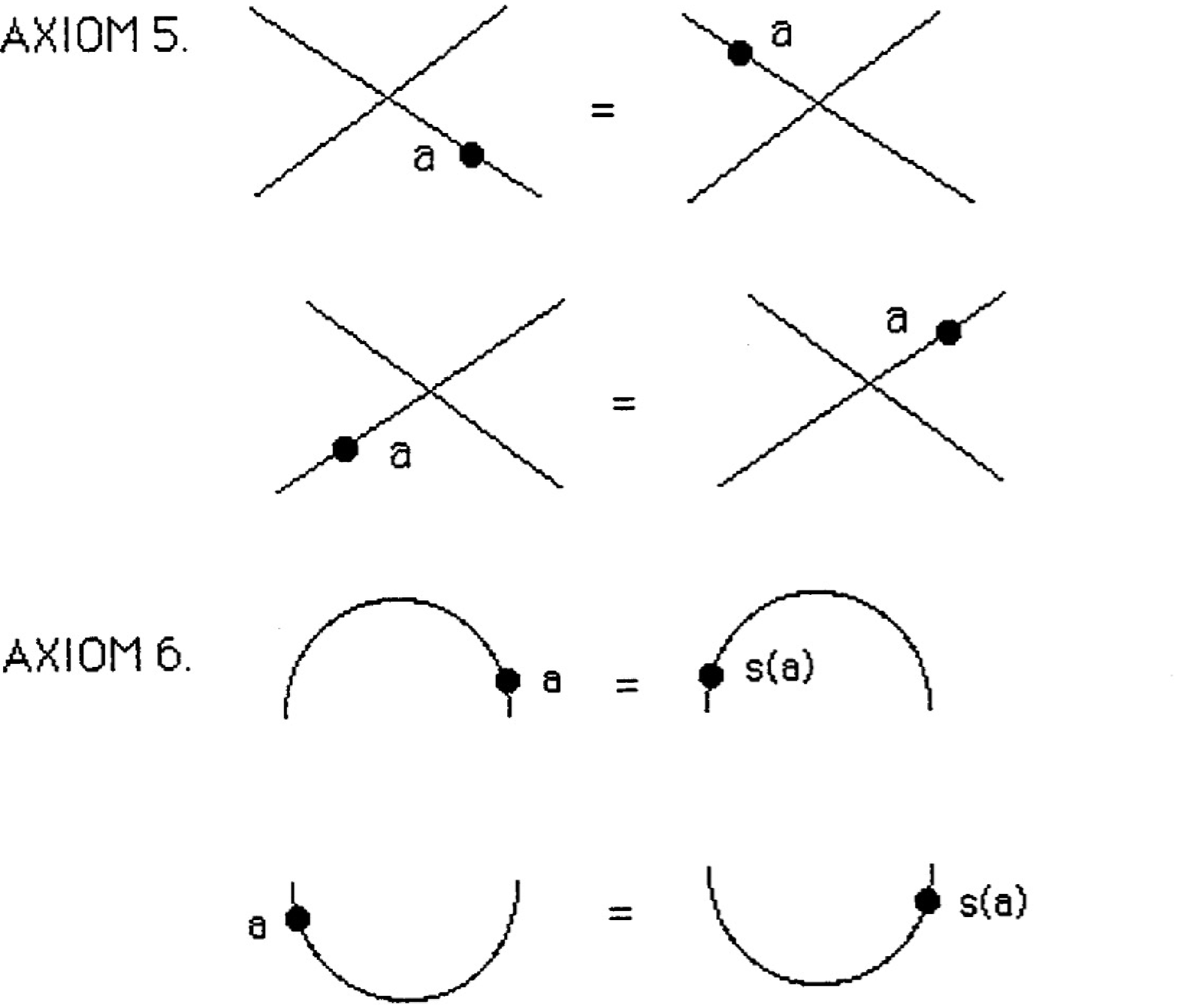}
     \end{tabular}
     \caption{\bf Permutation and Antipode}
     \label{F7}
\end{center}
\end{figure}

\subsection{The Category Arising From A Hopf Algebra} Let $A$ be an Hopf algebra. We shall define a category $Cat(A)$ associated with this Hopf algebra. $Cat(A)$ is a generalization of the immersion category $Flat$. In the case where $A$ is quasi-triangular, we shall define a functor $F:Tang \longrightarrow Cat(A)$. The invariants described in the later sections of this paper are consequences of the existence of this functor. \vspace{3mm}

The objects of $Cat(A)$ are identical to the objects of $Flat$, except that we take the abstract object $k$ and replace it by the ground ring of the Hopf algebra. Rather than make a separate notation for the distinction between $k$ as an abstract object and $k$ as the ground ring, we shall treat this contextually. Unless otherwise specified, a morphism from $k$ to $k$ is an abstract morphism, as in $Flat$. Each element of $A$ is taken to be a morphism from $V$ to $V$ and the composition $ab$ of elements $a$ and
$b$ in $A$ is simply their
product in $A$. Similarly, each element of the $n-fold$ tensor product of $A$ with itself is interpreted as a morphism from the $n$-fold tensor product of $V$ to itself. Each element $x$ of the ground ring $k$ is interpreted as a morphism from $k$ to $k$ in the same way, except that now we can take this morphism $x$ as right multiplication by $x$ if we wish. The generating morphisms of $Cat(A)$ consist in the morphisms of the tensor powers of $A$ (and $k$) together with the morphisms already available in $Flat$. The relations on the morphisms in $Flat$ still hold, and we add the following interrelationships with the morphisms from the Hopf algebra:
\vspace{3mm}

\noindent
5. For $a$ and $b$ in $A$ (viewed as morphisms of $V$ to $V$),
$$(a \otimes b)P = P(b \otimes a).$$ \\
6. For $a$ in $A$, let $s(a)$ denote the result of applying the antipode in $A$ to the element $a$. Then $$Cup(a \otimes 1_{V}) = Cup(1_{V} \otimes s(a))$$ and
$$(1_{V} \otimes a)Cap = (s(a) \otimes 1_{V})Cap.$$ \\ \vspace{3mm}

We have labeled these axioms $5$ and $6$, since the category Cat(A) already partakes of the flat tangle axioms $0$,$2$,$3$,$4$ with the caveat that $L$=$R$. Note that $5$ says that $P$ does act as a permutation on morphisms coming from $A$. See Figure 7 for an illustration of axioms 5
and 6.
\vspace{3mm}

There is one important addition to the structure of $Cat(A)$ that is not included in the tangle categories. In $Cat(A)$ we allow as morphisms formal sums of morphisms with coefficients in $k$. Thus in extending $\Delta(a) = \sum a_{1} \otimes a_{2}$ to become a morphism in $Cat(A)$, we take the sum of the individual morphisms in this
summation. A similar remark applies to the categorical interpretation of
identities such as $\sum \epsilon(a_{1}) a_{2} = a$. \vspace{3mm}

Axiom $6$ gives a direct relationship between the antipode in the Hopf algebra and the diagrammatic antipode that we described for the tangle category and flat tangle category. To see this relationship, we need to make a few remarks about the structure of the diagrams that represent morphisms in $Cat(A)$. A symbolic element $a$ of the Hopf algebra $A$ is diagrammed as a morphism in $Cat$ by taking a vertical line segment and labeling it with the letter $a$ next to a dot or "bead" drawn on the line segment. We will refer to the location of the bead on a larger diagram. Thus, in axiom $6$, $(1 \otimes a)Cap$ corresponds to a drawing of the $Cap$ with a bead labeled $a$ on its right hand side (below the maximum), while $(s(a) \otimes 1)Cap$ corresponds to a drawing of the cap with a bead labeled $s(a)$ on its left side. The equation $(1 \otimes a)Cap = (s(a) \otimes 1)Cap$ gives us permission to ``slide a bead counterclockwise around a maximum" while applying the antipode to it.
Similarly, the other equation of $6$ says that we can slide a bead counterclockwise around a minimum and apply the antipode to that bead.
Proofs by ``sliding" can often replace algebraic manipulations. For example, in Figure 8 we illustrate a proof of the following lemma (standard proof given below) of the expression of the antipode in terms of cups and caps. Figure 8 also illustrates the interpretation of sliding beads across a permutation $P$, and it gives a proof by sliding of the important formula $$s^{2}(a) = GaG^{-1}.$$ The square of the antipode (in $Cat(A)$) is given by conjugation with the formal grouplike $G$.
\vspace{3mm}

\noindent
{\bf Lemma.} Let $a$ be an element of $A$ viewed as a morphism from $V$ to $V$. Then $s(a) = (1_{V} \otimes Cup)(1_{V} \otimes a \otimes 1_{V})(Cap \otimes 1_{V})$. \\ \vspace{3mm}

\noindent
{\bf Proof.} $(1_{V} \otimes Cup)(1_{V} \otimes a \otimes 1_{V})(Cap \otimes 1_{V}) \\
=(1_{V} \otimes Cup)(1_{V} \otimes 1_{V} \otimes s(a))(Cap \otimes 1_{V})\\
=(1_{V} \otimes Cup) (Cap \otimes 1_{V})s(a)\\ =1_{V}s(a) = s(a).\\$
\vspace{3mm}

 \begin{figure}[htb]
     \begin{center}
     \begin{tabular}{c}
     \includegraphics[width=7cm]{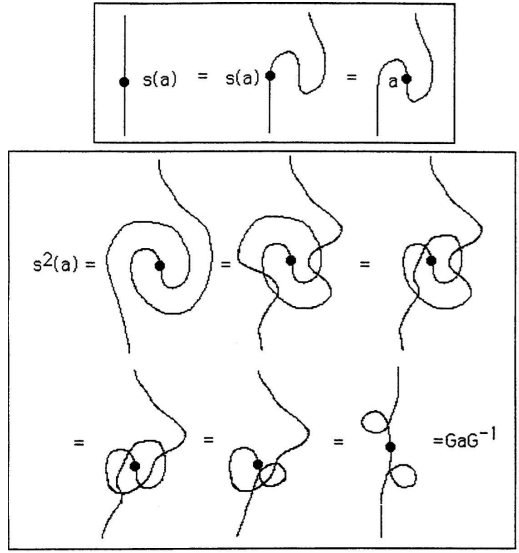}
     \end{tabular}
     \caption{\bf Proof by sliding.}
     \label{F8}
\end{center}
\end{figure}

Let $Cat_{1}(A)$ denote the set of all morphisms in $Cat(A)$ from $V$ to $V$ that are supported on a single strand flat tangle. In other words, a morphism in $Cat_{1}(A)$ consists of a single strand flat tangle that has been ``decorated" with beads from $A$ at various spots that are neither maxima, minima or crossings in the strand. It is clear from our axioms that we can slide all the beads through the tangle until they are all encountered first along a straight piece of the strand. Thus, if $T$ is this morphism, then $T$ is equivalent to a product $wt$ where $w$ is in $A$ and $t$ is a flat tangle free of elements of $A$. By our work with $Flat$, $t$ is equivalent to $G^{d}$ where $G$ is the flat curl $(1_{V} \otimes Cup)P(1_{V} \otimes Cap)$ and $d$ is the Whitney degree of $t$. Thus $T$ is equivalent to $wG^{d}$. This means that, except for closed loops, the morphisms in $Cat(A)$ are not much more than elements of the tensor powers of $A$ augmented by powers of the formal grouplike element $G$.
\vspace{3mm}

It may not be apparent at first sight that the element $w$ in $A$ that we obtained from $T$ is uniquely determined by the equivalence class of $T$ in $Cat(A)$. (It is clear from our previous remarks that the Whitney degree of $t$ is determined by the equivalence class of $T$.) In order to see this, we will give a definition of $w(T)$ that is dependent only on the decomposition of $T$ as a product of cups, caps, permutations and elements of tensor products of the algebra $A$. This is the same as saying that we will define $w(T)$ in terms of a given diagrammatic representation of $T$, since each diagrammatic representation is exactly a particular factorization of $T$ into elementary morphisms. \vspace{3mm}

The algorithm for computing $w(T)$ is as follows: Proceed upward along the strand of $T$, creating two data structures. The first data structure is an element in $A$ whose initial value is $1$. Let $w$ denote the generic form of this element of $A$. The second data structure is a power of the antipode of $A$ whose initial value is the identity mapping. Let $s^{i}$ denote the generic form of this power of the antipode. Now, whenever you encounter a "bead" $a$ on the strand, replace $w$ by $ws^{i}(a)$. Whenever you move around a maximum or a minimum in a clockwise direction, replace $s^{i}$ by $s^{i+1}$. Whenever you move around a maximum or a minimum in a counter-clockwise direction, replace $s^{i}$ by $s^{i- 1}$. The word $w(T)$ is the word obtained by going from the bottom of the strand to the top of the strand by this process. The final value of $i$ will be $d$, the Whitney degree of the tangle.
\vspace{3mm}

It is easy to see that $w(T)$ is invariant under all the replacements generated by the axiom for the category $Cat(A)$, and it is equally easy to see that $w(T)$ is exactly the element of $A$ that is obtained by sliding on a particular diagram. This shows that sliding is well-defined, and that the category $Cat(A)$ does not lose any information that is present in the Hopf algebra $A$. The Hopf algebra can be recovered from the morphisms of $Cat(A)$.
\vspace{3mm}

In the course of this discussion, we may have aroused the reader's appetite for the the structure of the closed loop morphisms in $Cat(A)$. This will be taken up in the next section, when we discuss trace and integral. We are now ready to construct a functor from $Tang$ to $Cat(A)$ when $A$ is quasi-triangular.
\vspace{3mm}

\subsection{The Functor F: $Tang \longrightarrow Cat(A)$ } Let $A$ be a quasi-triangular Hopf algebra as described in Section 1.
Let $\rho \in A \otimes A$ denote the Yang-Baxter element for $A$ and write $\rho$ symbolically in the form $\rho = \sum e \otimes e'$. We wish to define a functor from $Tang$ to $Cat(A)$. It suffices to define $F$ on the generating morphisms $R$,$L$,$Cup$ and $Cap$. We define
$$F(Cup) = Cup,$$
$$F(Cap) = Cap,$$
$$F(R) = P\rho = P (\sum e \otimes e')= \sum P (e \otimes e'),$$
$$F(L) = \rho^{-1} P = \sum (s(e) \otimes e') P = \sum (e \otimes s^{-1}(e')) P.$$
\vspace{3mm}

Diagrammatically, it is convenient to to picture $F(R)$ as a flat crossing with beads {\em above} the crossing labeled $e$ and $e'$ from left to right, with the summation indicated by the double appearance of the letter $e$. Similarly, $F(L)$ is depicted as a crossing with beads {\em below} the crossing and labeled $s(e)$ and $e'$. See Figure 9.
\vspace{3mm}

 \begin{figure}[htb]
     \begin{center}
     \begin{tabular}{c}
     \includegraphics[width=7cm]{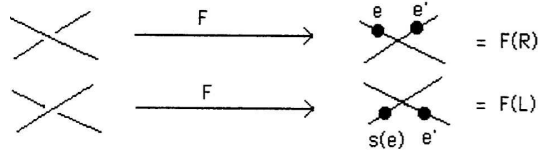}
     \end{tabular}
     \caption{\bf $F(R)$ and $F(L)$}
     \label{F9}
\end{center}
\end{figure}

Recall the axioms $0$,$2$,$3$ and $4$ for the tangle category.
\vspace{3mm}

\noindent
0. $$(Cup \otimes 1_{V})(1_{V} \otimes Cap) = 1_{V} \otimes 1_{V},$$ \\ 
$$(1_{V} \otimes Cup)(Cap \otimes 1_{V}) =1_{V}.$$ \\ 2. $$RL = LR = 1_{V} \otimes 1_{V}.$$\\
3. $$(R \otimes 1_{V})(1_{V} \otimes R)(R \otimes 1_{V}) = (1_{V} \otimes
R)(R \otimes 1_{V})(1_{V} \otimes R).$$\\ 4. $$(L \otimes 1_{V})( 1_{V}\otimes Cap) = (1_{V} \otimes R)(Cap \otimes 1_{V}).$$\\
\vspace{3mm}

In order for the functor $F$ to be well-defined, we must have the images of these relations satisfied in $Cat(A)$. For this, axiom $0$ follows at once since $Cup$ is taken to $Cup$ and $Cap$ to $Cap$. Axiom $2$ follows directly since we designed $F(L)$ as the inverse of $F(R)$. Axiom $3$ is exactly equivalent to the statement that $\rho$ satisfies the Yang-Baxter equation in $A$. Thus we are left to verify
\vspace{2mm}

\noindent
4. $$(F(L) \otimes 1_{V})( 1_{V} \otimes Cap) = (1_{V} \otimes F(R))(Cap \otimes 1_{V}).$$\\
\vspace{3mm}

This is verified by bead sliding in Figure 10. \vspace{3mm}

 \begin{figure}[htb]
     \begin{center}
     \begin{tabular}{c}
     \includegraphics[width=7cm]{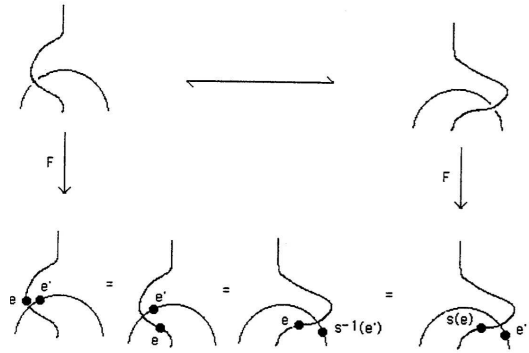}
     \end{tabular}
     \caption{\bf Switchback}
     \label{F10}
\end{center}
\end{figure}

\vspace{3mm}

Knowing that $F$ is a functor on the tangle category, means that we have implicitly defined many invariants of knots and links, since regularly isotopic tangles will have equivalent images under $F$. Before dealing with the intricacies of traces or of closed loops, we can state the following Theorem for 1-1 tangles, giving invariants with values in the Hopf algebra $A$.
\vspace{3mm}

\noindent
{\bf Theorem.} Let $T$ be a single-stranded, 1-1 tangle. That is , $T$ is a ``knot on a string". Then $F(T) = w[T]G^{d(T)}$ is a regular isotopy invariant of $T$. Here $w[T]=w(F(T))$ is the element of $A$ defined in the last section by concentrating all the algebra in $F(T)$ in the lower part of the tangle, and $d(T)$ is the Whitney degree of the plane curve underlying $T$. In fact $w[T] \in A$ is itself a regular isotopy invariant of $T$, as is $d(T)$.
\vspace{3mm}

\noindent
{\bf Proof.} This follows directly from the definition and well-definedness of the functor $F$ in conjunction with the discussion in the section on the category associated with a Hopf algebra.
\vspace{3mm}

\noindent
{\bf Example.} The curl {\em $v_{TOP}$} obtained from $G^{-1}$ by placing a crossing of type $L$ at its self-intersection maps, under $F$, to the ribbon element $v$ when $A$ is a ribbon Hopf algebra. $F(v_{TOP}) = v.$ The factorization of v into the product $G^{-1} \sum s(e')e$ is implicated by the slide convention for the antipode and the fact that $(s \otimes s) \rho = \rho.$ See Figure 11. \vspace{3mm}

 \begin{figure}[htb]
     \begin{center}
     \begin{tabular}{c}
     \includegraphics[width=7cm]{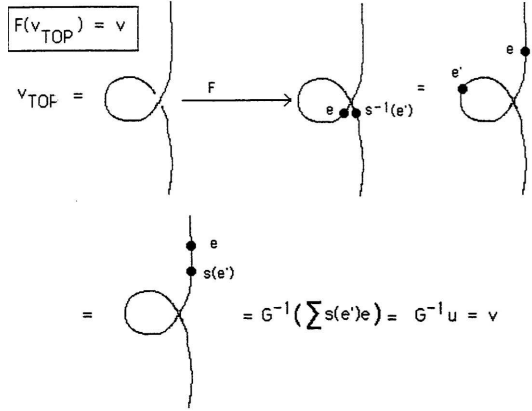}
     \end{tabular}
     \caption{\bf The Ribbon Element}
     \label{F11}
\end{center}
\end{figure}

\noindent
{\bf Remark.} When the identification $s(v_{TOP})=v_{TOP}$ is added to regular isotopy, the twists catalog only the framing, and the equivalence relation on the link diagrams is the same as ambient isotopy of framed links. See Figure 12. We call this equivalence
relation on link diagrams {\em ribbon equivalence.} Recall from Section 2 that a quasitriangular Hopf algebra is said to be a {\em ribbon Hopf algebra}
if there exists a grouplike element $G$ such that $v = G^{-1}u$ is in the center of $A$ and $s(u) = G^{-1}uG^{-1}$ where $u=\sum s(e')e.$ Note that algebraically, the
condition $s(v)=v$ is equivalent to the condition $s(u) = G^{-1}uG^{-1}.$
Thus Figures 11 and 12 show that a ribbon element is the exact counterpart
of ribbon equivalence under the functor $F.$ \vspace{3mm}

 \begin{figure}[htb]
     \begin{center}
     \begin{tabular}{c}
     \includegraphics[width=7cm]{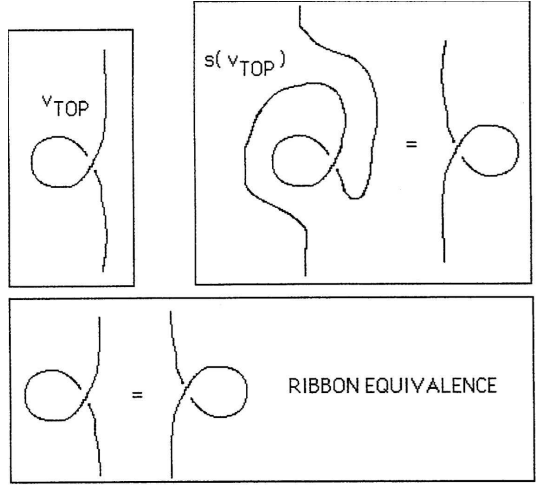}
     \end{tabular}
     \caption{\bf Ribbon Equivalence}
     \label{F12}
\end{center}
\end{figure}

\noindent
{\bf Remark.} In general, if $T$ is a single strand tangle, and $F(T)$ is
the corresponding element in the quasitriangular Hopf algebra $A$ determined by our correspondence, then
$F(\Delta (T)) = \Delta (F(T))$ where the first $\Delta$ is the diagrammatic coproduct and the second $\Delta$ is the algebraic coproduct. This fact follows from the axioms for a quasitriangular Hopf algebra. In particular, we use axiom 2 from Section 2:
$$\rho_{13} \rho_{12} = (1_{A} \otimes \Delta) \rho,$$ $$\rho_{13} \rho_{23} = (\Delta \otimes 1_{A})\rho,$$ and the fact (a consequence of the axioms) that $$\Delta(s(x)) = s(x_{2}) \otimes s(x_{1})$$ when $\Delta(x) = x_{1} \otimes x_{2}.$
The naturality of the coproduct with respect to the functor $F$ is then a
consequence of naturality with respect to the generating morphisms $Cup$,
$Cap$, $R$ and $L$ as illustrated in Figure 13. \vspace{3mm}

\begin{figure}[htb]
     \begin{center}
     \begin{tabular}{c}
     \includegraphics[width=7cm]{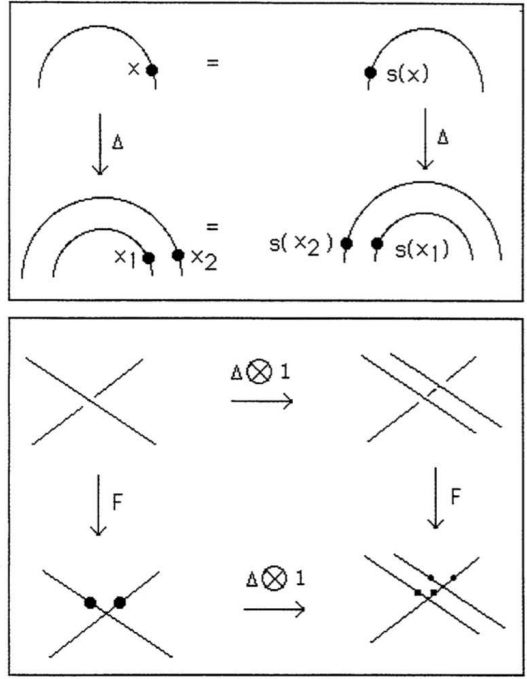}
     \end{tabular}
     \caption{\bf The Coproduct}
     \label{F13}
\end{center}
\end{figure}

\section{Diagrammatic Geometry and the Trace} An {\em augmented} Hopf algebra $A$ is a Hopf algebra that contains a grouplike element $G$ such that $s^{2}(x) = GxG^{-1}$ for all $x$ in $A$.
We can adjoin such an element to any given Hopf algebra (See \cite{Q} for
the details.). The result is an algebra in which the flat curl morphisms in
$Cat(A)$ can be identified with $G$ and $G^{-1}$. A ribbon Hopf algebra is
an augmented quasi-triangular Hopf algebra with the extra properties that
$s(u) = G^{-1}uG^{-1}$ and $v = G^{-1}u$ is in the center of $A$ where $u$ is the Drinfeld element described in Section 2.
It is useful to abstract the concept of augmented Hopf algebra when dealing
with the category, $Cat(A)$, of a given Hopf algebra $A$. We shall assume
throughout this section that the Hopf algebra $A$ is augmented. \vspace{3mm}

A function $tr:A \longrightarrow k$ from the Hopf algebra to the base ring $k$ is said to be a {\em trace} if
$tr(xy) = tr(yx)$ and $tr(s(x))= tr(x)$ for all $x$ and $y \in A.$ In this section we describe how a trace function on an augmented Hopf algebra yields an invariant, $TR(K)$, of regular isotopy of knots and links.
\vspace{3mm}

Let $ CMorph(A)$ denote the set of closed morphisms in $Cat(A)$ from $k$
to $k$. These morphisms are obtained from closed immersions of (collections) of circles that are decorated with elements of $A$ and arranged with respect to the vertical to form morphisms in $Cat(A)$ from $k$ to $k$. We would like to interpret such morphisms as actual mappings of
the base ring $k$ to itself. In $Cat(A)$ they are formal morphisms from $k$
to $k$ until further interpreted.
\vspace{3mm}

Define the {\em right and left circle morphisms} $$O_{R}:A \longrightarrow CMorph(A)$$ and $$O_{L}:A \longrightarrow CMorph(A)$$ by taking $O_{R}(a)$ for $a$ in $A$
to the the morphism obtained from a circle (simple composition of $cap$ and
$cup$) by placing $A$ on the right hand side of the circle. That is, 
$$O_{R}(a) = Cap(1_{V} \otimes a)Cup.$$ Similarly, $$O_{L}(a) = Cap(a \otimes 1_{V})Cup,$$ is the morphism that results from
placing $a$ on the left hand side of the circle. Note that by our axioms for $Cat(A)$ we have that
$O_{L}(a) = O_{R}(s^{-1}(a))$ (by sliding the $a$ bead around the circle).
\vspace{3mm}

Here are some fundamental identities about these left and right circle morphisms.

$$O_{L}(a) = O_{R}(s^{-1}(a)),$$
$$O_{R}(a) = O_{R}(s^{2}(a)),$$
$$O_{R}(a) = O_{L}(G^{-1}aG^{-1}) = O_{R}(s(a)G^{2}).$$ 

\noindent
These identities are proved diagrammatically in Figure 14. We will see that these identities are directly related to the formal structure of traces and integrals on a Hopf algebra. \vspace{3mm}

\begin{figure}[htb]
     \begin{center}
     \begin{tabular}{c}
     \includegraphics[width=7cm]{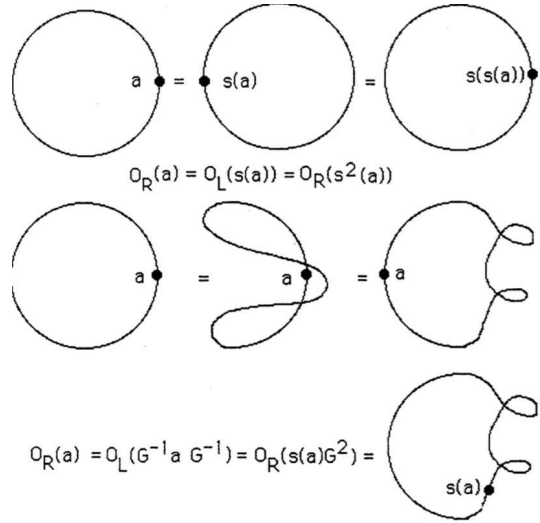}
     \end{tabular}
     \caption{\bf Circle Morphisms}
     \label{F14}
\end{center}
\end{figure}

\noindent
{\bf Slide Lemma.} If $W$ is any single-component closed morphism in $CMorph(A)$, then there is an element $a$ in $A$ such that $W= O_{R}(a)$. (Here equality denotes equality of morphisms in $Cat(A)$.)
\vspace{2mm}

\noindent
Writing the tensor product of closed morphisms as juxtaposition, any n-component closed morphism has the form $W = \sum O_{R}(a_{1})O_{R}(a_{2})...O_{R}(a_{n})$ for some elements $a_{1}, a_{2},...,a_{n}$ in $A$. In this expression the summation sign denotes the
possibility that there may be a summation over elements of $A$ that is shared among the components of the morphism. \vspace{3mm}

\noindent
{\bf Proof.} By sliding, concentrate the algebra on the single-component morphism W into $w \in A$ on a single segment of the immersion associated
with $W$. Now use the Whitney-Graustein Theorem to transform the immersion
associated with W to a circle decorate with a product of curls. Translate
the curls into a power of $G$ and amalgamate this with the algebra. The result is a circle decorated with an algebra element $a=wG^{k}$. Without
loss of generality, we can assume that $a$ is on the right side of the circle. This shows that
$W= O_{R}(a)$ as desired. The multi-component statement follows by the same
argument.
\vspace{3mm}

We can now point out a formal construction that has the properties of a trace. Define $\tau: A \longrightarrow CMorph(A)$ by the equation $$\tau(a)
= O_{R}(aG).$$
Note that another way to describe $\tau(a)$ is to say that the morphism $\tau(a)$ is obtained by inscribing $a$ on a curve with vanishing Whitney
degree (interpreting the $G$ in the definition of $\tau$ as a curl on the
circle). This means that if $a$ is slid all the way around this curve it
will return to its original position unchanged. \vspace{3mm}

\noindent
{\bf Lemma.} With $\tau$ described as above, $$\tau(ab) = \tau(ba)$$ for any $a$ and $b$ in $A.$ $$\tau(s(a))= \tau(a)$$ for any $a$ in $A$. \vspace{3mm}

\noindent
{\bf Proof.} $\tau(ab) = O_{R}(abG) = O_{R}(baG)$ by the remarks about sliding that precede the statement of this Lemma. Hence $\tau(ab) = O_{R}(baG) = \tau(ba).$
For the second part, $\tau(s(a)) = O_{R}(s(a)G) = O_{R}(s(s(a)G)G^{2}) =
O_{R}(s(G)s^{2}(a)G^{2}) =
O_{R}(G^{-1}s^{2}(a)G^{2}) = O_{R}(G^{-1}GaG^{-1}G^{2}) = O_{R}(aG) = \tau(a).$ (Alternatively, view Figure 15.) This completes the
proof.
\vspace{3mm}

\begin{figure}[htb]
     \begin{center}
     \begin{tabular}{c}
     \includegraphics[width=7cm]{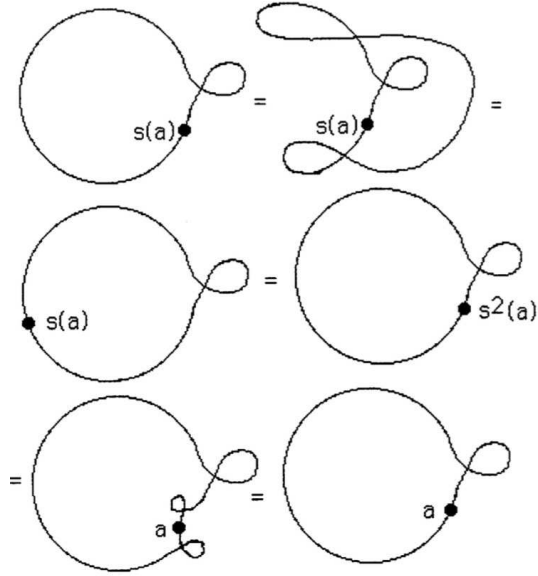}
     \end{tabular}
     \caption{\bf Tau Identities}
     \label{F15}
\end{center}
\end{figure}

\noindent
{\bf Remark.} Note that since $\tau(a) = O_{R}(aG)$, we have 

$$O_{R}(a) = \tau(aG^{-1})$$ and
$$O_{R}(a) = O_{L}(s^{-1}(a)) = \tau(s^{-1}(a)G).$$ 

\noindent
This remark tells us that any closed single-component morphism $W$ can be
expressed in terms of the formal trace $\tau$. For, by the Slide Lemma above, we can write $W = O_{R}(a) = \tau(aG^{-1})$. Recall from the proof
of the Slide Lemma that $a = wG^{k}$ where $w$ is the element of $A$ that
results from concentrating the algebra of $W$ to a single segment of $W$.
$G^{k}$ is the further concentration of curls at this segment that results
from applying the Whitney-Graustein Theorem to the immersion for $W$. {\em
If we then orient $W$ so that the selected segment has an upward arrow, then it is easy to see that $d = k-1$ is the Whitney degree of the immersion associated to $W$.}
Thus have the
\vspace{3mm}

\noindent
{\bf Evaluation Lemma.} For any closed single-component morphism $W$, $W =
\tau(aG^{d})$ where $a$ is the result of concentrating the algebra of $W$
on any segment of the immersion for $W$ and $d$ is the Whitney degree of this immersion, oriented so that the arrow is up at this segment. This formula for $W$ does not depend upon the choice of the segment where the concentration occurs. (The analogous result holds for a multi-component morphism. In this case, each closed loop can be dealt with separately as a
formal trace, and there will usually be a summation over products of these
traces.)
\vspace{3mm}

\noindent
{\bf Proof.} Most of the proof has already been given. Note that in concentrating the algebra, one may end up with $aG^{k}$ on the left hand side of the circle. The relationship $O_{L}(s^{-1}(x)) = O_{R}(x)$ then
shows that the same prescription works in this case. This completes the proof.
\vspace{3mm}

\noindent
{\bf Definition.} Two elements $a$ and $b$ of $A$, are said to be {\em slide equivalent} if one can be obtained from the other by either globally
applying the antipode or by rewriting the order of a product decomposition.

\vspace{2mm}

Thus $a$ and $s(a)$ are slide equivalent and $ab$ is slide equivalent to $ba$. (Note that $abcd$ is slide equivalent to $bcda$ but we make no assertion about the equivalence of $abcd$ and $acbd$.) Thus, by the properties of the formal trace $\tau$, we see that if $a$ is slide equivalent to $b$, then $\tau(a) = \tau(b)$. The next Lemma proves a converse to this statement.
\vspace{3mm}

\noindent
{\bf Recovery Lemma.} Let $a$ and $b$ be elements of $A$. Then $a$ and $b$
are slide equivalent if and only if
$\tau(a) = \tau(b)$.
\vspace{3mm}

\noindent
{\bf Proof.} Recall that $\tau(a) = O_{R}(aG^{-1})$. Thus $\tau(a)$ is obtained by decorating an immersion of Whitney degree zero with the algebra
element $a$. It is easy to see that slide equivalent elements in $A$ have
the same image under $\tau$. In fact, we have already formalized this fact
by proving that $\tau(xy) = \tau(yx)$ and that $\tau(s(x)) = \tau(x)$. Since
$\tau(a) = \tau(s(a))$ we can ``forget" about the applications of the antipode to $a$ as $a$ is moved across a maximum or a minimum in any curve
of total Whitney degree zero. That is, by our axioms that only way that $a$
can change in the course of equivalences to the morphism $O_{R}(aG^{-1})$
is by the application of the antipode to either all of $a$ or to some of its factors partially slid around the curve. Any time any factor is slid all the way around a curve it returns to its original value because the Whitney degree is zero. Then, since only slide equivalences are produced
by regular homotopies of this immersion to itself, $\tau(a) = \tau(b)$ implies that $a$ and $b$ are slide equivalent. This completes the proof.
\vspace{3mm}

\noindent
{\bf Remark.} Suppose that $tr:A \longrightarrow k$ is a trace function. That is, $tr$ is a linear function satisfying \vspace{3mm}

\noindent
1. $tr(xy)=tr(yx)$ and
\vspace{3mm}

\noindent
2. $tr(s(x))=tr(x).$
\vspace{3mm}

Then it follows from the Recovery Lemma that we may define tr on single-component closed morphisms $W$ by the formula $tr(\tau(a)) = tr(a)$ or
$tr(W) = tr(aG^{d})$ where $a$ and the Whitney degree $d$ are obtained by
concentrating the algebra on $W$ as in the Evaluation Lemma. This formula
extends to products and multicomponent closed morphisms in that obvious way. The upshot is that a trace function on the Hopf algebra gives rise to
an invariant of knots and links via our functor $F$. This is discussed in
the next section.
\vspace{3mm}

\noindent
{\bf Definition and Computation of TR(K).} Suppose that $tr:A \longrightarrow k$ is a trace function. In order to define an invariant of unoriented links, concentrate the algebra for each component of the link, and define $TR(K)$ to be the sum of the products of the evaluations of the individual components of the link. It follows from our previous discussion that this will be a regular isotopy
invariant of links.
\vspace{3mm}

\noindent
{\bf Discussion.} Note that invariants described according to the last definition can be regarded as computed either via the categorical decomposition of a morphism into cups, caps and crossings, or via the algebra concentration described in this section. In the categorical point
of view, a diagram with no free ends is a morphism from $k$ to $k$, where
$k$ is the ground ring of the Hopf algebra. It is useful to have both points of view available both for theory and for computation.

\section {Invariants of 3-manifolds}
The structure we have built so far can be used to construct invariants of 3-manifolds presented in terms of surgery on framed links. We sketch here our technique that simplifies an approach to 3-manifold invariants of Mark Hennings \cite{Hennings}. \vspace{3mm}

Recall that an element $\lambda$ of the dual algebra $A^{*}$ is said to be a {\em right integral } if $\lambda (x)1 = m(\lambda \otimes 1)(\Delta (x))$ for all $x$ in $A.$ For a unimodular \cite{LarsonSweedler} ,\cite{Radford-trace} finite dimensional ribbon Hopf algebra $A$ there is a right integral $\lambda$ satisfying the following properties for all x and y in A: \vspace{3mm}

\noindent
0) $\lambda$ is unique up to scalar multiplication when $k$ is a field.
\vspace{3mm}

\noindent
1) $\lambda (xy) = \lambda (s^{2}(y)x).$ \vspace{3mm}

\noindent
2) $\lambda (gx) = \lambda (s(x))$ where $g= G^{2}$, $G$ the special grouplike element for the ribbon element $v= G^{-1}u.$
\vspace{3mm}

Given the existence of this integral $\lambda$, define a functional $tr:A \longrightarrow k$ by the formula $$tr(x) = \lambda (Gx).$$
(It follows from the fact that $s^{2}(G)=G$ that $\lambda(Gx) = \lambda(xG).$)
\vspace{3mm}

It is then easy to prove the following theorem \cite{KandR}.
\vspace{3mm}

\noindent
{\bf Trace Theorem.} The function $tr$ defined as above satisfies \vspace{2mm}

\noindent
$tr(xy) = tr(yx)$ for all $x,y$ in$A.$
\vspace{2mm}

\noindent
and
\vspace{2mm}

\noindent
$tr(s(x)) = tr(x)$ for all $x$ in $A.$
\vspace{3mm}

The upshot of this theorem is that for a unimodular finite dimensional Hopf algebra there is a natural trace defined via the existent right integral. Remarkably, this trace is just designed to behave well with respect to handle sliding \cite{K-Hopf} , \cite{KandR}. Handle sliding is the basic transformation on framed links that leaves the corresponding 3-manifold obtained by framed surgery unchanged. See \cite{KirbyCalc}. This means that a suitably normalized version of this trace on framed links gives an invariant of
3-manifolds. For a link $K$, we let TR(K) denote the functional on links, as described in the previous section, defined via $tr$ as above.
\vspace{3mm}

To see how the condition on handle sliding and the property of being a right integral are related in our category, we refer the reader to Figure
16 where the basic form of handle sliding is illustrated and its algebraic
counterpart is shown. The algebraic counterpart arises when we concentrate
all the algebra in a given link component in one place on the diagram. The
component is then replaced by a circle and formally its evaluation is $O_{R}(x)$ for a suitable $x$ in the Hopf algebra. As the diagram shows, if
we let $\lambda(x) = O_{R}(x)$, then invariance under handle sliding is implicated by $\lambda$ being a right integral on the Hopf algebra. \vspace{3mm}

\begin{figure}[htb]
     \begin{center}
     \begin{tabular}{c}
     \includegraphics[width=7cm]{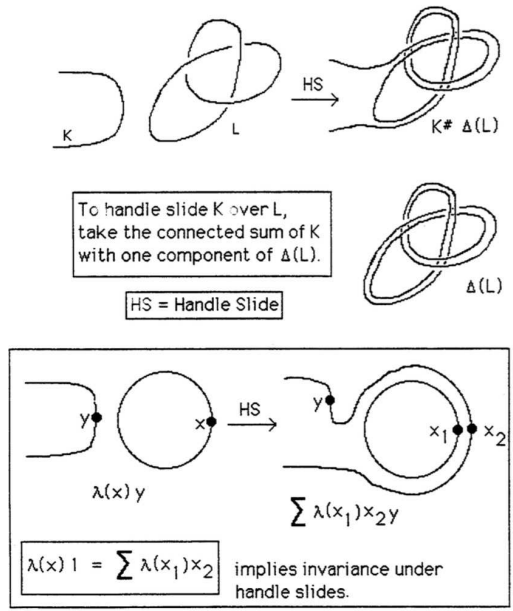}
     \end{tabular}
     \caption{\bf Handle Sliding and Right Integral}
     \label{F16}
\end{center}
\end{figure}

\vspace{3mm}

A proper normalization of $TR(K)$ gives an invariant of the 3-manifold obtained by framed surgery on $K.$ More precisely (assuming that $\lambda(v)$ and $\lambda(v^{-1})$ are non-zero), let

$$INV(K) = (\lambda(v) \lambda(v^{-1}))^{- c(K)/2}(\lambda(v) /
\lambda(v^{-1})^{- \sigma(K)/2}TR(K)$$

\noindent
where $c(K)$ denotes the number of components of K, and $\sigma(K)$ denotes the signature of the matrix of linking numbers of the components of $K$ (with framing numbers on the diagonal). Then $INV(K)$ is an invariant of the 3-manifold obtained by doing framed surgery on $K$ in the blackboard framing. This is our reconstruction of Hennings invariant \cite{Hennings} in an intrinsically unoriented context. \vspace{3mm}

\section{Centrality}
In the body of this paper we have described the functor $F$ which assigns
to each tangle a morphism in the category $Cat(A).$ If we choose a trace for each closed component of the tangle, Section 4 tells us how to map the
result of $F$ into the subcategory of $Cat(A)$ in which all components are
open (with a little care this second map, $F'(T),$ can be made into a functor). Finally, Section 3 describes how to map the result $F'(T)$ to an
element of $A^{\otimes n},$ where $n$ is the number of open strands. This
final quantity is a purely
algebraic invariant of the tangle and is the primary object of interest. In
particular, it assigns an element of $A$ to each 1-1 tangle. The range of
this invariant is not all of $A^{\otimes n}$ (or all of $A$ in the case of
1-1 tangles) however, as we show in this section. \vspace{3mm}

Notice that for the object $V^{\otimes n}$ in $Cat(A),$ every element $a$
of $A$ can be thought of as the morphism $\Delta^n(a)$ from the object to
itself (in particular, $a$ acts on $k$ as the morphism $\epsilon(a)$). Thus
$Hom(X,Y),$ the set of morphisms between two given objects $X$ and $Y,$ is
an $A$-bimodule, with $A$ acting on the left and right by left and right composition with the above morphism. Notice that the tensor product of morphisms is exactly the tensor product map of $A$-bimodules. We say that a
morphism
$x$ is {\it invariant\/} if it commutes with the action, so that $ax=xa$ for all $a$ in $A$ (the adjoint action of $A$ on these morphisms can easily
be defined in terms of this left and right action, and the definition of invariant says exactly that the morphism spans a trivial subrepresentation
in the adjoint representation, the usual definition of an invariant element). \vspace{3mm}

\noindent
{\bf Action Theorem.} The functor $F$ and the map $F'(T)$ from tangles to elements of $Cat(A)$ which have no closed components, as described above, take tan
gles to
invariant elements of $Cat(A).$
\vspace{2mm}

\noindent
{\bf Proof.} To show that $F(T)$ is invariant, since $F(T)$ is a composition of cups, caps and morphisms of the form $F(L)$, $F(R)$--images
of left and right crossings under the functor $F$--it suffices to show that
these morphisms commute with the action of $A$. However the statement that
$F(L)$ and $F(R)$ commute with this action is equivalent to the condition
$\rho \Delta = \Delta' \rho$ defining quasitriangularity. Finally, $Cap(x
\otimes y)g = \sum Cap(xg_{1} \otimes yg_{2}) = \sum Cap(xg_{1}s(g_{2}) \otimes y) = Cap(\epsilon (g) x \otimes y) = \epsilon(g) Cap(x \otimes y) =
g Cap(x \otimes y) .$ An identical argument applies to the $Cup$, and so we
conclude that $gF(T) = F(T)g.$ \vspace{2mm} 

Now in general Section 3 tells us that $F(T)$ can be written as $x \otimes
y,$ where $x$ is a product of closed components (i.e. a morphism in $Hom(k,k)$) and $y$ is an element of $Cat(A)$. Now since $A$ acts trivially
on $x,$ we have by the counit axiom that $a(x\otimes y)= x \otimes ay$ and
$(x \otimes y)a= x \otimes ay,$ and thus the invariance of $x \otimes y$ implies the invariance of $y$ (if $x$ is nonzero). But $F'(T)$ will send $x
\otimes y$ to $tr(x)y,$ where $tr$ is a product of one trace for each component, and thus to an invariant element. This concludes the proof. \vspace{3mm}

A 1-1 tangle $T$ is a tangle with a single input strand and a single output
strand. Then $F(T): V \longrightarrow V$, and by our axioms for sliding algebra around cups and caps, $F(T)$ is equivalent to the morphism corresponding to an algebra element of the form $ a(T) = wG^{d}$ where $G$
is the special grouplike element that we have discussed in the previous sections. (The element $a(T)$ is well-defined by sliding the algebra to the
bottom of the tangle and evaluating the closed loops in the tangle by a given functorial trace.)
\vspace{3mm}

\noindent
{\bf Centrality Theorem.} The algebra element $a(T)$ associated with a 1-1
tangle is in the center of the Hopf algebra $A$. \vspace{3mm}

\noindent
{\bf Proof.} By the {\em Action Theorem} we know that $gF(T) = F(T)g$ for
all $g$ in $A.$ Thus $ga(T) = gF(T)(1_{A})= F(T)g(1_{A}) = (F(T)(1_{A}))g =
a(T)g.$ Hence $a(T)$ is in the center of $A$. \vspace{3mm}

\noindent
{\bf Remark.} The argument that we have given to prove the centrality theorem may appear at first sight to be quite abstract. In fact, for each
example one can trace through the steps in the argument and produce a corresponding algebraic derivation of the commutativity that the theorem implies. Each stage in the categorical argument (applying an identity relating couinit and antipode or counit and coproduct, commuting a coproduct with an R-matrix) corresponds to an algebraic identity on the word(s) obtained by sliding all the algebra to the bottom of the diagram (using our sliding conventions in the category). The fact that each closed
component of the diagram is evaluated by a circle morphism means that the
evaluations obtained after sliding the algebra on a given closed curve (to
concentrate it in one segment) are independent of the location of the concentration.
\vspace{3mm}

A simple example (with no extra components) is illustrated in Figure 19. In this Figure the steps of the categorical proof are indicated, showing the successive locations of elements on the diagram. We can apply the functor $F$ to each of the diagrams in Figure 19, and then slide the algebra to the bottom of the diagram. Replacing the flat loops by $G$ or $G^{-1}$, we obtain an algebraic expression corresponding to each of the diagrams in the Figure. The algebraic expressions corresponding to each diagram in Figure 19 are listed below:

\begin{eqnarray}
aes^{-2}(e')G^{-1}\\
\epsilon(a_{1})a_{2}es^{-2}(e')G^{-1}\\
a_{2}es^{-1}(a_{11})s^{-2}(a_{12})s^{-2}(e')G^{-1}\\ ea_{12}s^{-1}(a_{11})s^{-2}(e')s^{-2}(a_{2})G^{-1}\\ e\epsilon(a_{1})s^{-2}(e')s^{-2}(a_{2})G^{-1}\\ es^{-2}(e')s^{-2}(a)G^{-1}&=& es^{-2}(e')G^{-1}a \end{eqnarray}

Each line in this list is a direct translation from the corresponding diagram in Figure 19, and each successive line is algebraically equal to its predecessor. Note how the diagrams supply the right powers of the antipode and other details that empower the resulting algebraic demonstration. These same principles apply to all examples, including the
case of extra components. We leave further examples as an exercise for the
reader, but recommend the case of a single line encircled once by an unknotted circle. The corresponding algebraic expression is $$f'e \lambda(fe')$$ where $\lambda$ stands for the circle morphism that is applied to the extra component. The exercise is to prove that this element
is in the center of the Hopf algebra by translating the categorical proof
to an algebraic proof. In the next section we shall take up this theme of
algebraic combinatorics again and show that our approach leads to other forms of algebraic proofs and to a deeper understanding of the nature of these central elements.
\vspace{3mm}

\noindent
{\bf Remark.} The reader should compare our treatment of centrality with \cite{Resh}. This centrality proof is remarkable in that it uses the structure of the category, $Cat(A)$, to prove an essentially algebraic fact about $A$. Of course it is in this category that $1-1$ tangles correspond to morphisms that are products of cups,caps and crossings. It is the structure of these building blocks that insures that tangles yield central elements in the algebra. To see the power of this argument, the reader should note that it proves that the ribbon element $G^{-1}u$ is in the center of the algebra, and that this in turn proves (via the categorical representation of the square of the antipode as conjugation by the grouplike element $G$) that the square of the antipode is represented by conjugation by the Drinfeld element $u$. This proof that $s^{2}(x) = uxu^{-1}$ is quite different from the direct algebraic proof.
\vspace{3mm}

\noindent
{\bf Remark.} A natural question is whether there are elements of the center
of the Hopf algebra which are not in the range of the functor $F$ applied to
1-1 tangles. In the cases of most interest, including the finite-dimensional
quantum groups at roots of unity, the answer is no. These and many other
cases where there is a Hennings-type invariant have the following property,
referred to by Hennings as unimodularity: The morphism $\rho P \rho P$ of
$Cat(A)$ corresponding
to a full twist of two strands, corresponds via $F$ to an element of $A \otimes
A,$ and therefore can be viewed as a map from $A^*$ to $A$. Unimodularity
means this map
is nondegenerate. In particular, if $f$ is $tr( \cdot G)$ for some trace
$tr,$ this map sends $f$ to the value of $F$ on the 1-1 tangle consisting of a
single vertical strand encircled by a closed loop evaluated with $tr.$ By the
preceding theorem, the resulting element of $A$ is in fact in the center.
Thus we have a 1-1 map from traces to the image of $F.$ It is shown in \cite{Radford-Kauf} that the space of traces on $A$ has the same dimension
as the
center of $A,$ and thus that this map is actually onto the center. Thus we
conclude that when the Hopf algebra is unimodular, $F$ maps 1-1 tangles onto the
center. A subtler question is whether the image of {\bf single strand} 1-1
tangles generates the center.
\vspace{3mm}

\section{Centrality, Algebra and Combinatorics} 

It is puzzling that centrality is proved so smoothly using the categorical
structure when it appears to be quite intricate at the level of pure algebra. The purpose of this section is to show how the computations appear
at the algebra level and how this level is related to the combinatorics of
link diagrams. In particular, we will finish the section with a new proof
of centrality that confirms a longstanding conjecture of the first two authors of the paper. We conjectured that a certain algebraic method of moving an element $a$ across a sum of words $W$ would always result in a verification that $aW = Wa$ for $W$ in the image of our functor. As we shall
see, this method is directly related to the content of the categorical proof of centrality, and it has interesting combinatorial properties of its own.
\vspace{3mm}

Let $\rho = e \otimes e'$ denote the Yang-Baxter element for a quasi-triangular Hopf algebra $A.$ Let $s$ denote the antipode of $A$. In calculations below we follow the modified summation convention for indices. Thus $\Delta(a) = \Sigma a_{1} \otimes a_{2}$ will be simply denoted by $a_{1} \otimes a_{2}.$ Other examples of this usage are the formulas $a_{1} \epsilon (a_{2}) = a$ and $s(a_{1})a_{2} = \epsilon(a).$
\vspace{3mm}

Note the following calculation:

\begin{eqnarray}
ae \otimes e' &=& a_{2}e \otimes \epsilon(a_{1})e' \nonumber \\ &=& a_{3}e \otimes s(a_{1})a_{2}e' \nonumber \\ &=& ea_{2} \otimes s(a_{1})e'a_{3}. \nonumber \\ \end{eqnarray}

Thus we have the identity
$$ae \otimes e' = ea_{2} \otimes s(a_{1})e'a_{3}.$$ 

\begin{figure}[htb]
     \begin{center}
     \begin{tabular}{c}
     \includegraphics[width=7cm]{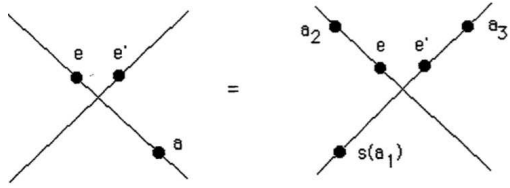}
     \end{tabular}
     \caption{\bf Bead Push Identity}
     \label{F17}
\end{center}
\end{figure}

In Figure 17, we illustrate the diagram corresponding to this identity. This diagram suggests that we could ``see" how a given element of the Hopf
algebra is in the center by pushing beads on the diagram. For example, consider the following calculation with
$$X=es^{-2}(e')G^{-1}.$$

\begin{eqnarray}
aX &=& aes^{-2}(e')G^{-1} \nonumber \\
&=& ea_{2}s^{-2}(s(a_{1})e'a_{3})G^{-1} \nonumber \\ &=& ea_{2}s^{-1}(a_{1})s^{-2}(e')s^{-2}(a_{3})G^{-1} \nonumber \\ &=& e\epsilon(a_{1}) s^{-2}(e')s^{-2}(a_{2})G^{-1} \nonumber \\ &=& e\epsilon(a_{1})s^{-2}(e')G^{-1}a_{2}GG^{-1} \nonumber \\ &=& es^{-2}(e')G^{-1} \epsilon(a_{1})a_{2} \nonumber \\ &=& es^{-2}(e')G^{-1} a \nonumber \\
&=& Xa. \nonumber \\
\end{eqnarray}

This is exactly the sort of intricate algebraic argument that our categorical proof of centrality seems to avoid. Now view Figure 18. In Figure 18, we show that the above algebraic proof of centrality has an exact diagrammatic counterpart via the identity from Figure 17. \vspace{3mm}

\begin{figure}[htb]
     \begin{center}
     \begin{tabular}{c}
     \includegraphics[width=7cm]{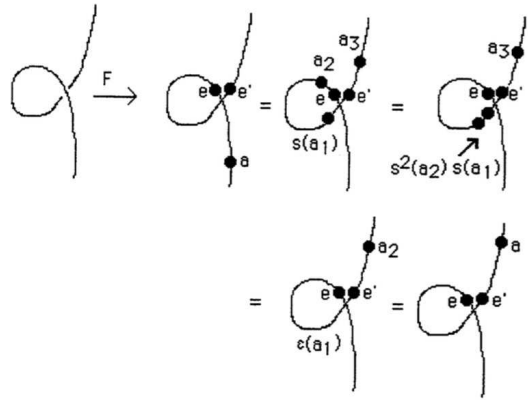}
     \end{tabular}
     \caption{\bf Centrality by Pushing Beads}
     \label{F18}
\end{center}
\end{figure}

Now view Figure 19. In this Figure we have illustrated the centrality of
the same element as in Figure 18, but the diagrammatic proof follows the pattern of the Centrality Theorem of the last section. We see from this example that the proof of the Centrality Theorem actually does provide a sequence of algebraic steps that gives a specific proof of centrality for
any given element of the Hopf algebra that is an image of a $1-1$ tangle $T$ under the functor $F$. The algebraic proof that is so constructed is
guided by the diagram of the $1-1$ tangle as it is arranged with respect to
a vertical direction. The steps in the algebraic proof parallel the movement of the element $a$ across the sequence of morphisms into which $F(T)$ is decomposed.

\begin{figure}[htb]
     \begin{center}
     \begin{tabular}{c}
     \includegraphics[width=7cm]{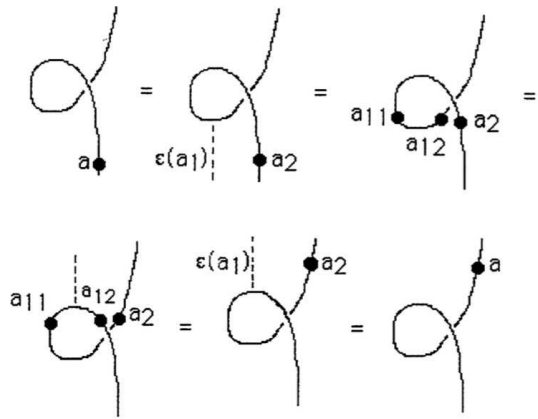}
     \end{tabular}
     \caption{\bf Algebraic Centrality via the Category}
     \label{F19}
\end{center}
\end{figure}

We leave it to the reader to translate the diagrams of Figure 19 into an
algebraic demonstration. The point is that once a given diagram is chosen,
then the steps of moving elements across the elementary morphisms are exactly specified and each step in this process yields an algebraic step that can be verified by the usual means. \vspace{3mm}

We now show that the method corresponding to Figures 17 and 18 will always
work to provide proofs of centrality.
In this method, we generalize the formula in Figure 17 to all the different
cases of a bead at one of the legs of a crossing. It is easy to verify that the resulting patterns do not depend upon the crossing type. The exact result is given in the next Lemma, whose proof we omit. See Figure
20 for the pattern of diagrammatic bead slides. In this figure the crossings are indicated with a dark vertex that can be
{\em either} an undercrossing or an overcrossing. Note that if the bead has an antipode applied to it, then the order of indices will shift from clockwise around the crossing to anticlockwise around the crossing (or vice-versa). With the help of the diagrams in Figure 20 one can experiment
on link diagrams and produce proofs of centrality by direct bead pushing as
we did in Figure 17. View Figures 21 for an example of this procedure. \vspace{3mm}

\begin{figure}[htb]
     \begin{center}
     \begin{tabular}{c}
     \includegraphics[width=7cm]{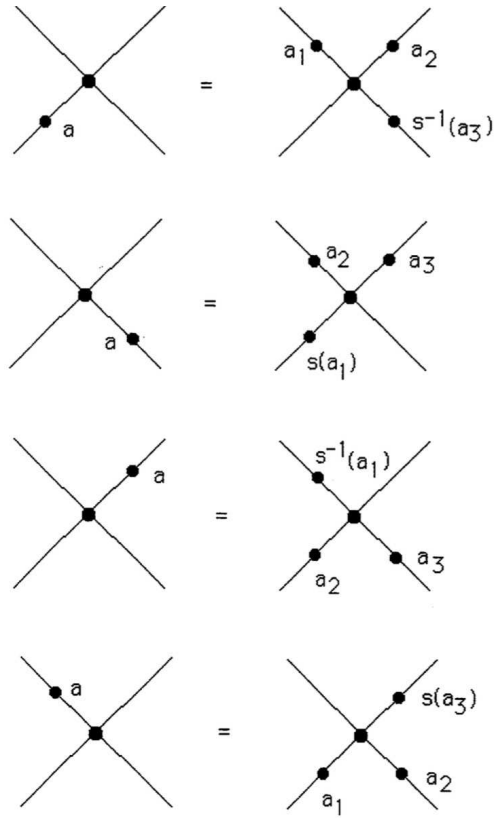}
     \end{tabular}
     \caption{\bf Bead Push Patterns}
     \label{F20}
\end{center}
\end{figure}

\noindent
{\bf Bead Push Lemma.} Let $\rho = e \otimes e'$ denote the Yang-Baxter
element for a quasi-triangular Hopf algebra $A.$ Let $s$ denote the antipode of $A$. (In the following formulas $\Sigma a_{1} \otimes a_{2}$
will be denoted by
$a_{1} \otimes a_{2}.$)
\vspace{3mm}

\noindent
1. $ae \otimes e' = ea_{2} \otimes s(a_{1})e'a_{3}$ \vspace{2mm}

\noindent
2. $e \otimes ae' = s^{-1}(a_{3})ea_{1} \otimes e'a_{2}$ \vspace{2mm}

\noindent
3. $ea \otimes e' = a_{2}e \otimes a_{1}e's(a_{3})$ \vspace{2mm}

\noindent
4. $e \otimes e'a = a_{3}es^{-1}(a_{1}) \otimes a_{2}e'$ \vspace{2mm}

\noindent
5. $as(e) \otimes e' = s(e)a_{2} \otimes s^{-1}(a_{3})e'a_{1}$ \vspace{2mm}

\noindent
6. $s(e) \otimes ae' = s(a_{1})s(e)a_{3} \otimes e'a_{2}$ \vspace{2mm}

\noindent
7. $s(e)a \otimes e' = a_{2}s(e) \otimes a_{3}e's^{-1}(a_{3})$ \vspace{2mm}

\noindent
8. $s(e) \otimes e'a = a_{1}s(e)s(a_{3}) \otimes a_{2}e'$ \vspace{2mm}

\noindent
{\bf Proof.} Omitted.
\vspace{3mm}

\begin{figure}[htb]
     \begin{center}
     \begin{tabular}{c}
     \includegraphics[width=7cm]{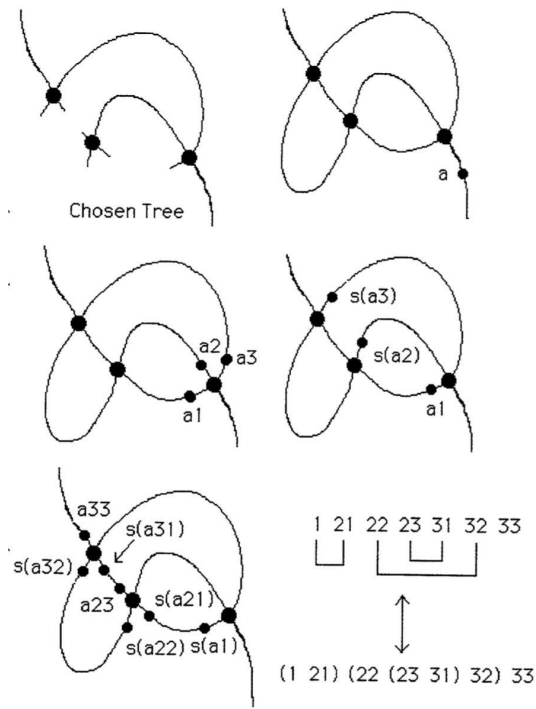}
     \end{tabular}
     \caption{\bf Trefoil Tree Push}
     \label{F21}
\end{center}
\end{figure}

In Figure 21 we illustrate one way to verify centrality for a trefoil tangle via bead pushing (here illustrated without crossing choices since the patterns of bead pushing do not depend upon the choice of crossing). In this procedure, we choose a connected tree that is obtained from the tangle $T$ for the trefoil by cutting midpoints of an appropriate subset of
the edges of the projected flat tangle. Once this tree is chosen, there is
a unique way to push a bead (labeled $a$) from the lower leg of the tangle
to all the branches of the tree including the upper tangle leg. Except for
the upper and lower legs of the tangle all the twigs of the tree are paired
by the cutting arcs.
We see in this example that each pair of paired twigs gives rise to a ``cancellation" of the form $s(a_{1})a_{2} = \epsilon(a)$ and ``reconstruction" in the form $\epsilon(a_{1})a_{2} = a$. These cancellations occur in sequence via the lexicographic ordering corresponding to the particular splitting into coproducts that is dictated
by the tree. Thus in the case of Figure 21 we have the ordering $$(1,21),(22,(23,31),32),33$$

\noindent
The parentheses indicate the pairings. Thus $1$ and $21$ are paired; $23$
and $31$ are paired; $22$ and $32$ are paired. If $x$ and $y$ are paired
then we set a left parenthesis before $x$ and a right parenthesis after $y$
in the lexicographic ordering. In the Figure 21 we have replaced the parentheses by connecting arcs. Now note that the parentheses so obtained
are nested in the classical fashion of well-formed parentheses. A cancellation of $23,31$ then leads to a cancellation of $22,32$ and there
is a parallel cancellation of $1,21.$ In the end we are left with $a$ reconstructed at the top edge of the tangle and hence a proof of centrality
for this particular tangle.
\vspace{3mm}

We now wish to show that the example in Figure 21 is quite general and that
this procedure will work on any $1-1$ tangle. In order to accomplish this
end it must be shown that
\vspace{2mm}

\noindent
1. For any choice of tree in a tangle $T$, each pair of paired beads will
have powers of the antipode applied to them that differ by one when they are moved into the same vertical sector of a common edge (See Figure 18 for an example.).
\vspace{2mm}

\noindent
2. The lexicographic order of coproducts combined with the pairings gives
rise to a well-formed structure of parentheses (so that the cancellation and reconstruction can proceed).
\vspace{3mm}

Condition 1 is proved by noting that paired beads are part of a circuit in
the tangle, as illustrated in Figure 22. The (Figure 20) rules for bead pushing make it easy to see that the total exponent for going around this
circuit is the same as its Whitney degree(namely plus or minus one) just as we discussed in earlier sections. The key observation (available directly from Figure 20) is that the pattern of application of the antipode
in the bead push is exactly the same as if the crossing were smoothed horizontally or vertically and the bead was pushed (possibly across a maximum or a minimum) according to the rules of our category). \vspace{3mm}

\begin{figure}[htb]
     \begin{center}
     \begin{tabular}{c}
     \includegraphics[width=7cm]{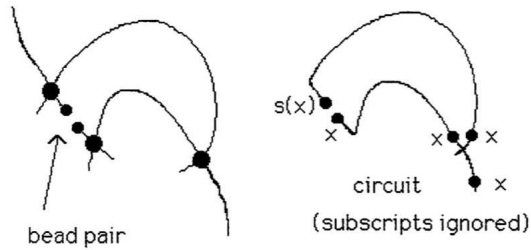}
     \end{tabular}
     \caption{\bf Condition 1}
     \label{F22}
\end{center}
\end{figure}

\begin{figure}[htb]
     \begin{center}
     \begin{tabular}{c}
     \includegraphics[width=7cm]{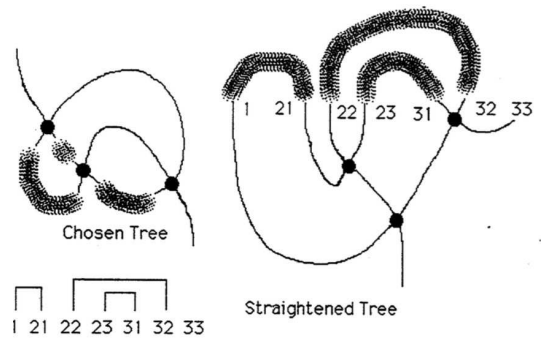}
     \end{tabular}
     \caption{\bf Condition 2}
     \label{F23}
\end{center}
\end{figure}

Condition 2 is proved by first noting (via the Bead Push Lemma) that the lexicographic ordering and pairings for a given choice of tangle diagram
and tree is the {\em same} as the lexicographic ordering and pairings in
the new diagram obtained by straightening the tree by a planar isotopy so
that each push moves upwards (with respect to the chosen vertical direction) and there are no maxima in the arcs of the tree. In this case
each vertex in the tree gives a clockwise ordering as dictated by the Bead
Push Lemma and Figure 20. In Figure 23 we illustrate this isotopy and the
resulting ordering. The reader should compare this ordering with the ordering in Figure 21. Once the tree is straightened, it is clear that the
well-formedness of the parenthesis structure corresponds to the fact that
the arcs connecting paired beads (each such arc is now a maximum with respect to the vertical) do not intersect one another. The well-formed parentheses are a direct consequence of the planarity of the graph of the
flat tangle.
\vspace{3mm}

This completes the proof that bead pushing into the twigs of a tree will always prove centrality for $1-1$ tangles. The reader should note that these arguments apply to tangles with multiple components just so long as
the circle morphisms for closed components take values in the ground ring
$k$ (See Section 4.). We have spent the effort to relate bead pushing with
centrality because the result is intriguing and because it shows quite clearly the relationship between centrality in the Hopf algebra and the combinatorics of plane graphs and trees. \vspace{3mm}

It remains to be seen if these methods can be inverted to give a complete
characterisation of central elements in quasi-triangular Hopf algebras (i.e. elements that correspond to the image of 1-1 tangles under our functor from tangles to the category associated with a Hopf algebra). \vspace{3mm}

It is worth remarking that the result we have described is invertible at the formal level in the following sense: If we know how to thread (that is prove $aW = Wa$) an arbitrary element $a$ of the Hopf algebra through a given sum of words $W$ where $W$ as a summation has the structure
of a sum over repeated Yang-Baxter elements, and if this threading involves
only applications of the bead push identities and the counit and antipode
identities (as in our discussion), then the threading will, of its own accord, produce a tree with a lexicographic ordering of the branches and a
legal parenthetical association of these branches. This data gives an embedding of the tree and parenthesis arcs into the plane. That planar embedding specifies a link diagram whose image is $W$ under our functor.
\vspace{3mm}

\noindent
{\bf A Quantum Remark.} One can interpret the ``beads" moving on the tangle
diagrams as "particles" that scatter through a ``spin network" that corresponds to the given tangle. In this interpretation, each element $a$
of the Hopf algebra (seen as a morphism in $Cat(A)$) is a quantum particle
traveling forward in time. The application of the antipode, $s(a)$ is interpreted as ``$a$ traveling backwards in time." The identity $\epsilon(a_{1})a_{2} = a$ is interpreted as the emission or
absorption by $a$ of a ``virtual photon". Note that when the virtual photon
$\epsilon(a_{1})$ is present, then $a$ is in a mixed state connoted by the
summation over $a_{1}$ and $a_{2}$.
The identity $s(a_{1})a_{2} = \epsilon(a)$ corresponds to the creation or
annihilation of a particle and an antiparticle. (Note that grammatically it
is not ``a particle and an antiparticle" but rather a ``particle/antiparticle
mixed state". Then our result on centrality is interpreted by saying that
the $1-1$ tangle is a ``self-energy diagram" in which the particle undergoes
a virtual interaction that returns it to its original state. 

\section{Applications to Virtual Knot Theory}
The categorical methods in this paper can be applied to virtual knot theory. The flat crossings of the category of a Hopf algebra have the same formal properties as the virtual crossings
in rotational virtual knot theory. See \cite{RV} for a detailed explanation of rotational virtual knot theory and the construction of the generalization of the methods of this paper for that category.
In defining a functor on the virtual tangle category one takes the virtual crossing in a tangle diagram to a flat crossing
in the Hopf algebra category. Another possible application to virtual knots for the three manifold invariants in this paper can be investigated via the virtual three manifolds, Kirby calculus and invariants of three manifolds as constructed in \cite{DK}. This application to invariants of virtual three-manifolds will be the subject of a further paper.\\

\end{document}